\pdfoutput=1
\documentclass[pdflatex,sn-mathphys-num]{sn-jnl} 

\usepackage{graphicx}%
\usepackage{multirow}%
\usepackage{amsmath,amssymb,amsfonts}%
\usepackage{amsthm}%
\usepackage{mathrsfs}%
\usepackage[title]{appendix}%
\usepackage{xcolor}%
\usepackage{textcomp}%
\usepackage{manyfoot}%
\usepackage{booktabs}%
\usepackage{subcaption}%
\usepackage{algorithm}%
\usepackage{algorithmicx}%
\usepackage{algpseudocode}%
\usepackage{listings}%
\usepackage{longtable}
\usepackage{soul}
\usepackage{bm}
\usepackage{dblfloatfix}%
\usepackage{placeins}%
\usepackage{afterpage}%
\usepackage{amsfonts}
\usepackage{dsfont}
\DeclareUnicodeCharacter{02BC}{'}

\theoremstyle{thmstyleone}%
\theoremstyle{thmstyletwo}%

\theoremstyle{thmstylethree}%

\newcommand{\danique}[1]{\textcolor{black}{#1}}

\newcommand{\set}[1]{\mathcal{#1}}

\begin{document}
\title[Rapid drone-based wildfire detection at a fraction of current prevention spending]{Rapid drone-based wildfire detection at a fraction of current prevention spending}


\author*[1]{\fnm{Romain} \sur{Puech}}\email{puech@mit.edu}

\author[2]{\fnm{Danique} \sur{de Moor}}\email{demoor@mit.edu}

\author[3]{\fnm{Ana} \sur{Trišović}}\email{ana\_tris@mit.edu}


\author[2]{\fnm{Dimitris} \sur{Bertsimas}}\email{dbertsim@mit.edu}

\affil*[1]{\orgdiv{Operations Research Center}, \orgname{Massachusetts Institute of Technology}, \orgaddress{\street{77 Massachusetts Avenue}, \city{Cambridge}, \postcode{02139}, \state{MA}, \country{USA}}}

\affil[2]{\orgdiv{Sloan School of Management}, \orgname{Massachusetts Institute of Technology}, \orgaddress{\street{77 Massachusetts Avenue}, \city{Cambridge}, \postcode{02139}, \state{MA}, \country{USA}}}

\affil[3]{\orgdiv{Computer Science \& Artificial Intelligence Laboratory}, \orgname{Massachusetts Institute of Technology}, \orgaddress{\street{77 Massachusetts Avenue}, \city{Cambridge}, \postcode{02139}, \state{MA}, \country{USA}}}

\abstract{\danique{Early wildfire detection is critical to prevent small ignitions from escalating into large-scale disasters, yet current monitoring systems lack a quantitative framework for allocating detection infrastructure at scale. We jointly optimize the placement of monitoring infrastructure and the routing of autonomous drones under realistic operational constraints to quantify the investment required for rapid, large-scale wildfire detection. Evaluated out-of-sample on 3,693 California ignitions from 2021–2024, an optimized drone network operating at a \$100 million five-year budget detects 97.3\% of fires, including 74\% within the first hour. Amortized over five years, that budget is about \$20 million per year, roughly 5\% of California's annual wildfire-prevention expenditure. Under current technology costs, drone-based monitoring is substantially more cost-effective than static ground sensors. Detection is governed primarily by spatial coverage, while routing strategy mainly determines detection speed.}}

\maketitle




\section*{Introduction}

\noindent Wildfires and associated forest loss have intensified in recent years, causing widespread ecological and economic damage~\citep{curtis2018classifying, tyukavina2022global, hu2025coexposure}. In 2025 alone, wildfires burned approximately 390 million hectares of land worldwide~\citep{GAR2025}. Between 2014 and 2023, wildfires caused an estimated \$106 billion in direct damage losses globally, not accounting for indirect losses linked to ecological degradation, biodiversity loss, and human displacement. A disproportionate share of these losses occurred in the United States~\citep{GAR2025}, especially California. 

Recent wildfire disasters highlight the increasing severity of these events. Nearly half of the 200 most damaging wildfires since 1980 have occurred in the past decade~\citep{science25}. The 2019–2020 ``Black Summer'' fires in Australia burned nearly 19 million hectares and affected an estimated three billion animals~\citep{nhra2023blacksummer,ward2020biodiversity}. In the United States, the 2018 Camp Fire caused 85 fatalities and near-complete destruction of the town of Paradise \cite{calfire2025deadliest, hamideh2022wildfire}; across that year, California's wildfires produced an estimated \$148.5 billion economic footprint, including capital losses, health costs, and supply-chain disruption \cite{wang2021economic}. More recently, the January 2025 Los Angeles wildfires destroyed over 16,000 structures and displaced tens of thousands of residents \cite{laedc2025wildfires}. These disasters coincide with increasingly extreme fire weather linked to climate change and the expansion of human settlements in the wildland–urban interface, which together amplify wildfire risk~\citep{science25}. Looking forward, empirical models project that even under about 1.5 °C warming, extra-tropical regions will experience larger and more intense wildfires, and current suppression policies are projected to fail across much of the world~\cite{haas2026wildfires}.


Effective wildfire management depends critically on how quickly new ignitions are detected. Fires identified during their initial stages are substantially easier and less costly to contain, whereas even short delays can allow them to spread rapidly, increasing suppression costs, ecological damage, and risks to nearby communities. Investments in early detection can therefore generate substantial societal and economic benefits. Realizing these benefits, however, requires monitoring systems capable of rapidly detecting ignitions across large and heterogeneous landscapes.

Existing wildfire monitoring systems rely on a combination of ground-based sensors and satellites, each providing valuable but inherently incomplete coverage. Ground-based monitoring systems provide continuous observations but are costly to deploy at scale. Satellite observations offer broad geographic reach but may be constrained by revisit frequency, cloud cover, smoke, or spatial resolution during the earliest stages of a wildfire \citep{WildfireDetectionBouguettaya}. Unmanned aerial vehicles (UAVs) offer a complementary approach by combining broad spatial coverage with the flexibility to continuously adapt surveillance to evolving wildfire risk~\citep{puech2026wfdronebench}. Recent advances in autonomous flight, sensing technologies, and onboard artificial intelligence have further increased the potential of UAV-based monitoring~\citep{WildfireDetectionDLKim,WildfireDetectionDLKumar24,WildfireDetectionDLYang}. The remaining task is to determine how monitoring infrastructure should be deployed and how aerial patrols should be coordinated so that new ignitions can be found quickly under realistic operational constraints.
\danique{Previous work has studied sensor placement~\citep{DroneroutingDeLaFuente,fire7070254}, drone routing~\citep{WildfireMonitoringBailonruiz22,DroneroutingPostGhamry17,DroneroutingXu,Tzoumas2023,DroneroutingJemmali,Demir2024,Pordal2025}, integrated infrastructure planning~\citep{Liu2025,puech2026wfdronebench}, and wildfire suppression optimization~\cite{boussioux2026predictiveprescriptiveaioptimizing}, yet quantitative analyses of proactive wildfire monitoring at statewide scale remain limited, particularly those jointly considering infrastructure deployment, operational constraints, and drone routing strategies.}

To address this gap, using California as a representative case study, we evaluate three routing strategies and quantify how infrastructure deployment, operational constraints, and investment influence wildfire detection performance. Although California serves as the motivating application, the proposed methodology is general and applicable to other regions.

Our results demonstrate that rapid, near-complete wildfire detection is achievable within hours under realistic operational and budget constraints. An optimized drone-based monitoring network operating at a \$100 million budget detects 97.3\% of benchmark wildfire ignitions, including 74\% within the first hour of ignition and all but one detected fire identified within four hours.
Under current cost assumptions, drone-based monitoring substantially outperforms ground sensors, providing greater spatial coverage per unit cost. More broadly, these results reveal a fundamental system insight: wildfire detection performance is primarily determined by the spatial reach of the monitoring network, while routing strategies mainly influence detection speed once fires fall within range. This distinction highlights the central role of infrastructure design in enabling effective large-scale wildfire monitoring.



\danique{Our contributions are fivefold. First, we show that rapid, near-complete wildfire detection across California is achievable with an investment well below current wildfire-prevention spending. Second, we characterize the budget--detection frontier, revealing a threshold relationship between investment and spatial coverage. Third, we identify the technology-cost regimes in which drone, ground-sensor, and hybrid monitoring systems are optimal. Fourth, we develop two optimization-based routing strategies and, by comparing them with a baseline, show that infrastructure placement determines detection, whereas routing primarily determines detection speed. Finally, we demonstrate that an optimized \$100 million drone network exceeds an optimistic upper bound on the coverage of California's deployed ALERTCalifornia camera system.}

\section*{Results}
\subsection*{Study design}

\begin{figure*}[!b]
    \centering

        \includegraphics[width=\textwidth,height=0.68\textheight,keepaspectratio]{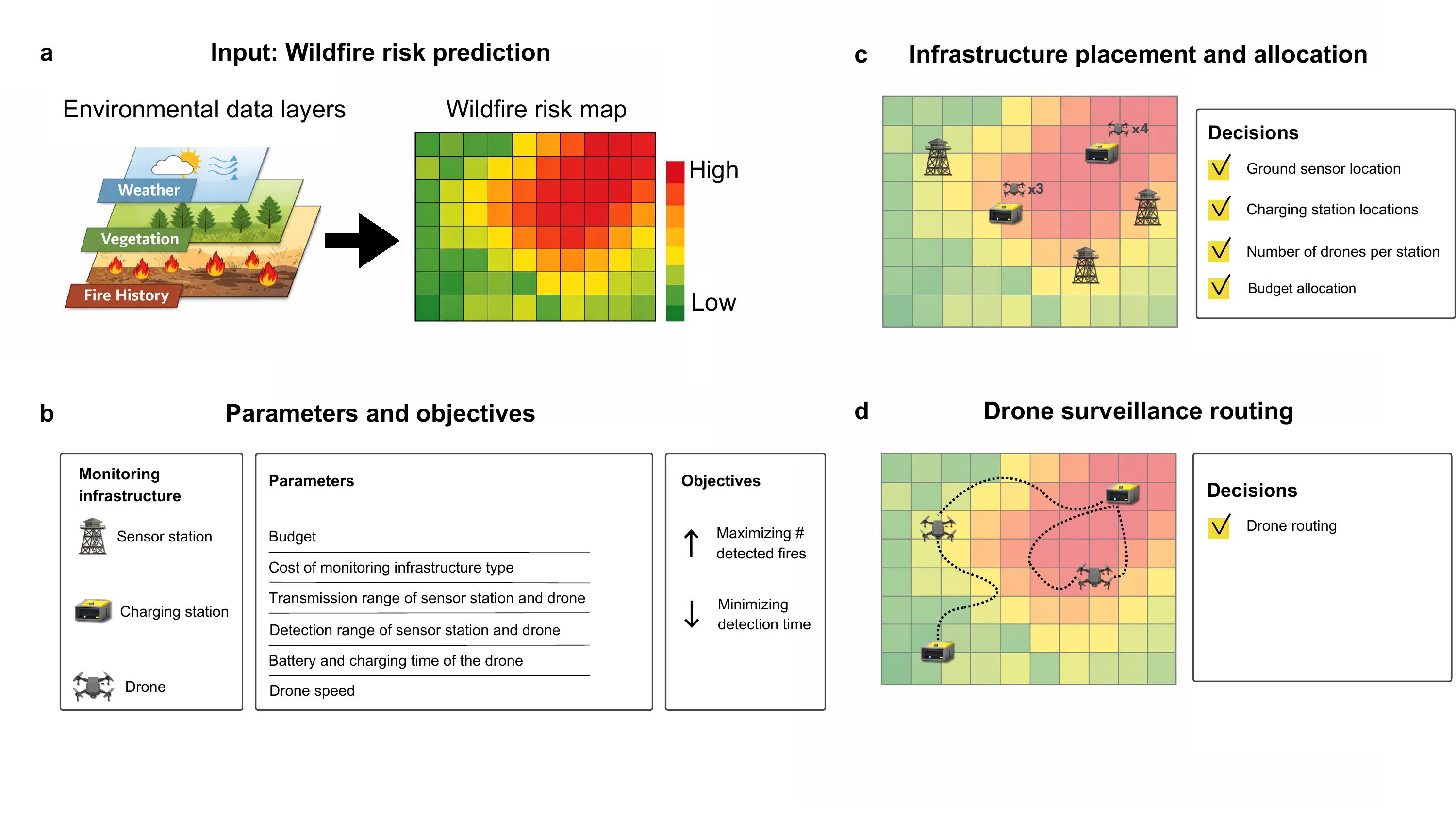}


    \caption{\textbf{AI-driven monitoring framework for proactive wildfire detection.} (a) Existing machine-learning models generate spatial wildfire risk maps from environmental data layers, including weather conditions, vegetation characteristics, and historical fire observations. (b) Key operational parameters and optimization objectives define the decision-making framework. (c) Based on the predicted wildfire risk, the monitoring system determines where to deploy sensor  and charging stations, as well as the number of drones allocated to each charging station to form a distributed monitoring network, using budget and infrastructure cost parameters. (d) Drones autonomously patrol the landscape following optimized surveillance routes that account for operational constraints including battery capacity, flight time limits, speed, communication range, and charging requirements.}
    \label{fig:framework}
\end{figure*}

To investigate the infrastructure requirements for rapid wildfire detection at statewide scale, we conduct a case study based on California, combining wildfire risk prediction, infrastructure optimization, and autonomous drone surveillance in an integrated monitoring system. The simulation environment and the underlying monitoring architecture follow~\citep{puech2026wfdronebench}, which benchmarks algorithms on synthetic instances. 


The monitoring system consists of ground sensors, charging stations, and drones. Infrastructure placement and drone routing are jointly optimized using predicted wildfire risk while accounting for budget limitations, battery endurance, communication range, charging requirements, and sensing capabilities (Fig.~\ref{fig:framework}).

\captionsetup[subfigure]{
    labelformat=simple,
    labelsep=space,
    font=bf,
    skip=1pt
}
\renewcommand\thesubfigure{\textbf{\alph{subfigure}}}

\subsection*{Near-complete wildfire detection within an hour}
\begin{figure*}[!t]
    \centering
        \includegraphics[width=\linewidth,height=0.62\textheight,keepaspectratio]{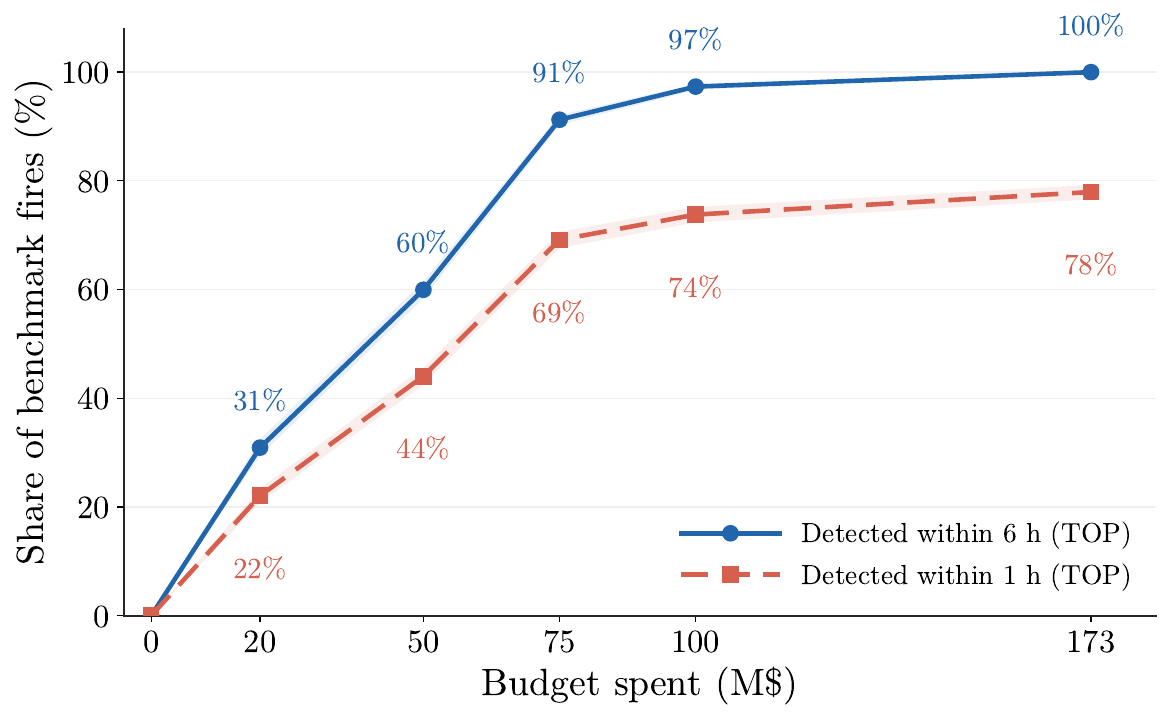}
        \caption{\textbf{Detection rate and speed as a function of monitoring investment.} Percentage of California wildfire ignitions across 2021--2024 ($n = 3{,}693$ fires, pooled) detected within six hours of ignition (blue) and within the first hour ($\Delta t = 0$, orange) for the team orienteering (TOP) routing strategy, as a function of realized hardware spend. Detection rates are computed over the full pooled fire population across all four years rather than by averaging per-year rates. Five budget levels are shown: \$20M, \$50M, \$75M, \$100M, and \$500M (from which only \$172.6M is spent). Shaded bands are 95\% Wilson confidence intervals.}
        \label{fig:frontier}
\end{figure*}

\danique{\noindent With a \$100 million budget, our optimized drone monitoring network detects 97.3\% (95\% CI: 96.7--97.8\%) of California wildfire ignitions across 2021--2024, with 74\% detected within the first hour of ignition ($\Delta t = 0$; Fig.~\ref{fig:frontier}, Table~\ref{tab:detection}). All but one detected fire was identified within four hours. Amortized over the five-year cost horizon, this \$100 million budget is equivalent to about \$20 million per year, or roughly 5\% of California's 2024--25 wildfire resource management and prevention expenditure~\cite{lao_wildfire_faq_2025}. Near-complete wildfire detection across California is therefore mainly a question of how much coverage to buy, given that we treat a fire as found once a drone is close enough.}

\begin{figure*}[!tb]
    \centering \includegraphics[width=\textwidth,height=0.66\textheight,keepaspectratio]{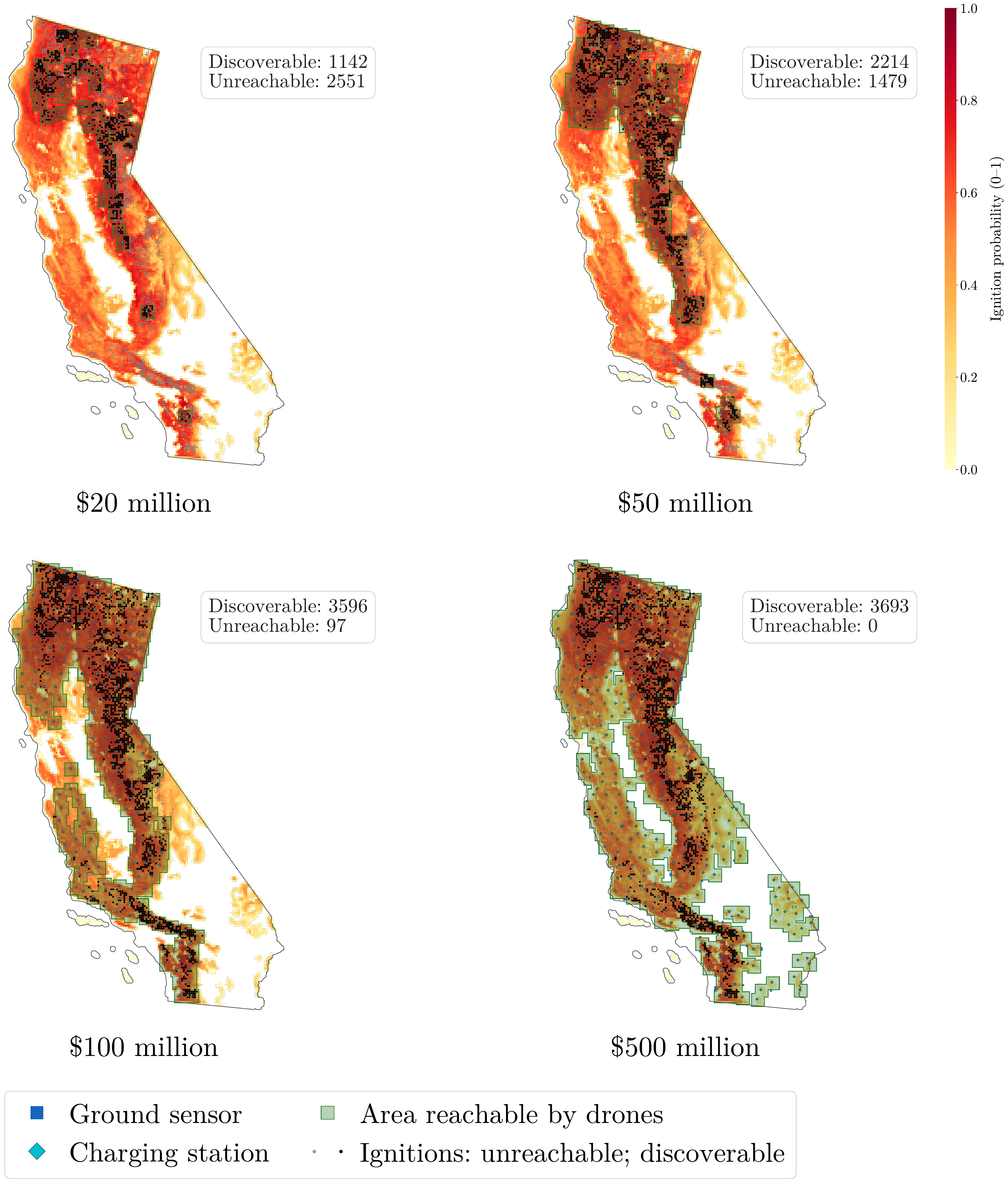}
    \caption{\textbf{Optimized deployment of monitoring infrastructure and fire discoverability across budget scenarios in California.} Fires are labeled as discoverable if the drones can reach them given their battery starting from any charging station. Discoverability is computed over all 3{,}693 historical California ignitions across 2021--2024. At \$20~million, 40~stations and 280~drones are deployed, concentrated in Northern and Central California; 31\% of historical fires are discoverable. At \$50~million, 100~stations and 700~drones expand coverage to 60.0\% of benchmark fires. At \$100~million, 200~stations and 1{,}400~drones span the majority of the state's burnable landscape; 97\% of historical fires are discoverable. At the \$500~million budget, 339~charging stations, 31~ground sensors, and 2{,}373~drones are deployed at a realized hardware cost of \$172.6~million; all historical ignitions are discoverable.}
    \label{fig:placement}
\end{figure*}

\danique{The detection frontier shows a pronounced nonlinear relationship between investment and detection (Fig.~\ref{fig:frontier}). Detection increases from 30.9\% at \$20 million to 60.0\% at \$50 million and 91.2\% at \$75 million, before plateauing above approximately \$100 million as additional investment yields diminishing returns. This threshold behavior reflects expanding spatial coverage (Fig.~\ref{fig:placement}): the fraction of historical ignitions outside the network's reach decreases from 69\% at \$20 million to only 3\% at \$100 million (3,596 of 3,693 fires discoverable). Once a fire is reachable, detection is nearly guaranteed, with rates exceeding 95\% for the best-performing strategy across all budget levels (Table~\ref{tab:detection}). At the \$500 million budget, where every historical ignition is discoverable, the sole missed fire results from an inconsistency between datasets: the USFS ignition record places the fire on a lake, where our risk map assigns zero ignition probability.}

\begin{table}[!t]
\caption{Detection performance across routing strategies and deployment budgets,
  pooled over California wildfire years 2021\textendash2024 (\(n = 3{,}693\) fires).
  Rates are computed over the full pooled fire population.
  \emph{Det.\ rate}: overall detection rate with 95\% Wilson score interval.
  \emph{Within 1\,h}: share of all fires detected in the first simulation step (\(\Delta t = 0\)).}
\label{tab:detection}
\begin{tabular*}{\textwidth}{@{\extracolsep{\fill}}llrr@{}}
\toprule
Budget & Strategy & Det.\ rate (\%) & Within 1\,h (\%) \\
\midrule
\$20M & TOP     & 30.9 [29.5, 32.4] & 22.1 \\
 & MaxCov         & 29.5 [28.1, 31.0] & 17.2 \\
 & LinearMinTime  & 25.9 [24.5, 27.3] & 18.4 \\
\addlinespace
\$50M & TOP     & 60.0 [58.4, 61.5] & 44.1 \\
 & MaxCov         & 57.0 [55.4, 58.6] & 32.7 \\
 & LinearMinTime  & 46.1 [44.5, 47.8] & 33.0 \\
\addlinespace
\$75M & TOP     & 91.2 [90.3, 92.1] & 69.1 \\
 & MaxCov         & 87.1 [86.0, 88.1] & 51.1 \\
 & LinearMinTime  & 79.9 [78.6, 81.2] & 56.9 \\
\addlinespace
\$100M & TOP     & 97.3 [96.7, 97.8] & 73.8 \\
 & MaxCov         & 92.7 [91.8, 93.5] & 53.9 \\
 & LinearMinTime  & 84.5 [83.3, 85.6] & 61.8 \\
\addlinespace
\$500M & TOP     & 100.0 [99.8, 100.0] & 77.9 \\
 & MaxCov         & 96.6 [96.0, 97.1] & 60.2 \\
 & LinearMinTime  & 84.3 [83.1, 85.4] & 64.9 \\
\bottomrule
\end{tabular*}
\end{table}

\subsection*{Drones dominate ground sensors economically}
\noindent \danique{Across all budget levels below \$500 million, the placement optimizer selects no ground sensors, allocating the entire infrastructure budget to charging stations and drones. Under our baseline cost assumptions (Methods), a charging station equipped with seven drones costs \$0.50 million and can scan approximately 1{,}225~km$^2$ per battery cycle, at roughly one-tenth the cost per unit area of an equivalent ground-sensor deployment. These assumptions are consistent with publicly reported deployment costs (Methods), including per-camera installation costs of \$15,000--\$100,000 for the Oregon camera network~\cite{KLCC2026OregonCameras} and annual operating costs of up to \$36,000 per station in utility-scale deployments~\cite{MauiNews2024HawaiianElectric}. Ground sensors appear only at the \$500 million budget to close the small residual coverage gaps remaining after statewide drone coverage is achieved (Fig.~\ref{fig:placement}).}

\begin{figure*}[!tb]
    \centering
    \includegraphics[width=\linewidth,height=0.60\textheight,keepaspectratio]{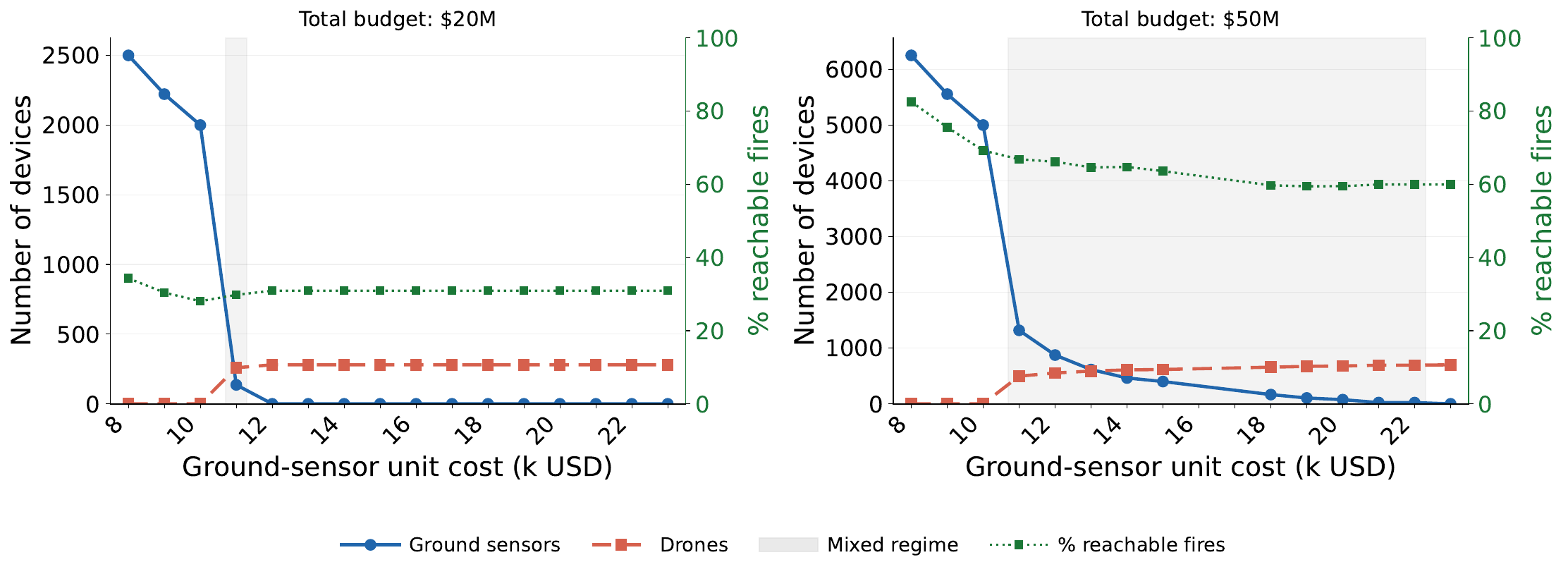}
    \caption{\textbf{Cost sensitivity of the optimal infrastructure mix to ground-sensor unit cost.} Optimal device counts as the assumed ground-sensor unit cost is varied, at total budgets of \$20~million (left) and \$50~million (right), with charging-station and drone unit costs held fixed at the values in Table~\ref{tab:benchmark_params}. \textbf{Left axis:} number of ground sensors (blue, solid) and drones (orange, dashed) in the optimal placement. \textbf{Right axis} (green, dotted): share of California wildfire ignitions across 2021--2024 reachable by the resulting network ($n = 3{,}693$ fires). The shaded band marks the \emph{mixed regime}, where ground sensors, charging stations, and drones coexist in the optimum: it is narrow at \$20~million (around USD~11k per sensor, switching to a drones-only portfolio by USD~12k) and substantially wider at \$50~million (roughly USD~11k--22k, drones-only from USD~23k).}
    \label{fig:breakeven_costsensitivity}
\end{figure*}

\paragraph{Cost sensitivity.}
To evaluate the sensitivity of our results to sensor costs, we varied only the unit cost of ground sensors while holding drone and charging-station costs fixed at the values in Table \ref{tab:benchmark_params}, and re-solved the placement optimization model as in our main experiments. At total budgets of \$20 million and \$50 million, we observe three qualitative regimes in the optimal infrastructure mix as the ground sensor unit cost varies (Fig.~\ref{fig:breakeven_costsensitivity}). When static sensors are very inexpensive, the budget is devoted almost entirely to dense ground coverage and no charging stations are deployed. At intermediate unit costs, ground sensors, charging stations, and drones coexist: aerial infrastructure extends reachable surveillance while static units remain where they are still comparatively cost-effective. When static units become sufficiently expensive, the optimal portfolio uses only charging stations and drones. The price interval over which mixed deployments arise is narrow at \$20 million (illustrated between approximately \$10{,}000 and \$12{,}000 per ground sensor in our scan) but substantially wider at \$50 million, where mixed portfolios persist from roughly \$11{,}000 to about \$22{,}000 per sensor. Our benchmark parameterization of \$100{,}000 per ground sensor (Table~\ref{tab:benchmark_params}) lies far above this coexistence band, which is consistent with the drone-heavy infrastructure choices reported above. If that \$15{,}000 lower bound were used as the five-year unit cost, the \$50 million optimum would be mixed; once operating costs are included, as in our total-cost-of-ownership figure, even the low end of the reported capital range lies above the coexistence band. Because the optimal networks below \$500 million contain no ground sensors, detection performance depends on aerial hardware costs and budget only through their ratio: uniformly scaling charging-station and drone costs by a factor $\alpha$ is equivalent to scaling the budget by $1/\alpha$. The detection frontier in Fig.~\ref{fig:frontier} therefore doubles as a sensitivity analysis for aerial hardware costs: even if our unit costs are optimistic by a factor of two, near-complete detection is achievable at ~\$200 million (about \$40 million per year amortized), still about 0.03\% of the 2018 damages~\cite{wang2021economic}.

\subsection*{Routing strategy primarily affects detection speed}
\noindent We evaluate three drone routing strategies: team orienteering (TOP), which maximizes the cumulative risk of cells visited per battery cycle; maximum coverage (MaxCov), which prioritizes surveillance of high-risk areas; and detection-time minimization (LinearMinTime), which concentrates patrols on the highest-risk cells to reduce detection delay (details in Methods). Across all budget levels, TOP achieves the highest overall detection rates, with the difference most pronounced at intermediate budgets. At \$100 million, TOP detects 97.3\% of benchmark fires versus 92.7\% for MaxCov and 84.5\% for LinearMinTime (Table~\ref{tab:detection}).

The key structural finding is that, for the best-performing routing strategies, placement governs detection rate and routing governs detection speed. Among fires within drone reach, detection is nearly guaranteed by the two best strategies; what varies is how quickly the fire is found. From Table \ref{tab:detection}, we observe that at \$100 million, 74\% of fires are detected by TOP within the first hour ($\Delta t = 0$), compared with 54\% for MaxCov (Mann-Whitney $U$ test on $\Delta t$ values among jointly detected fires, $p < 0.001$). LinearMinTime, despite being explicitly designed to minimize detection delay, detects fewer fires overall because its patrol concentration leaves peripheral locations undervisited. Placement is therefore the binding constraint: extending the network's geographic footprint matters more than how drones are routed within it.

\subsection*{Optimized drones outperform an optimistic upper bound on the ALERTCalifornia system's coverage}
\begin{figure*}[!t]
    \centering
    \includegraphics[width=\textwidth,height=0.62\textheight,keepaspectratio]{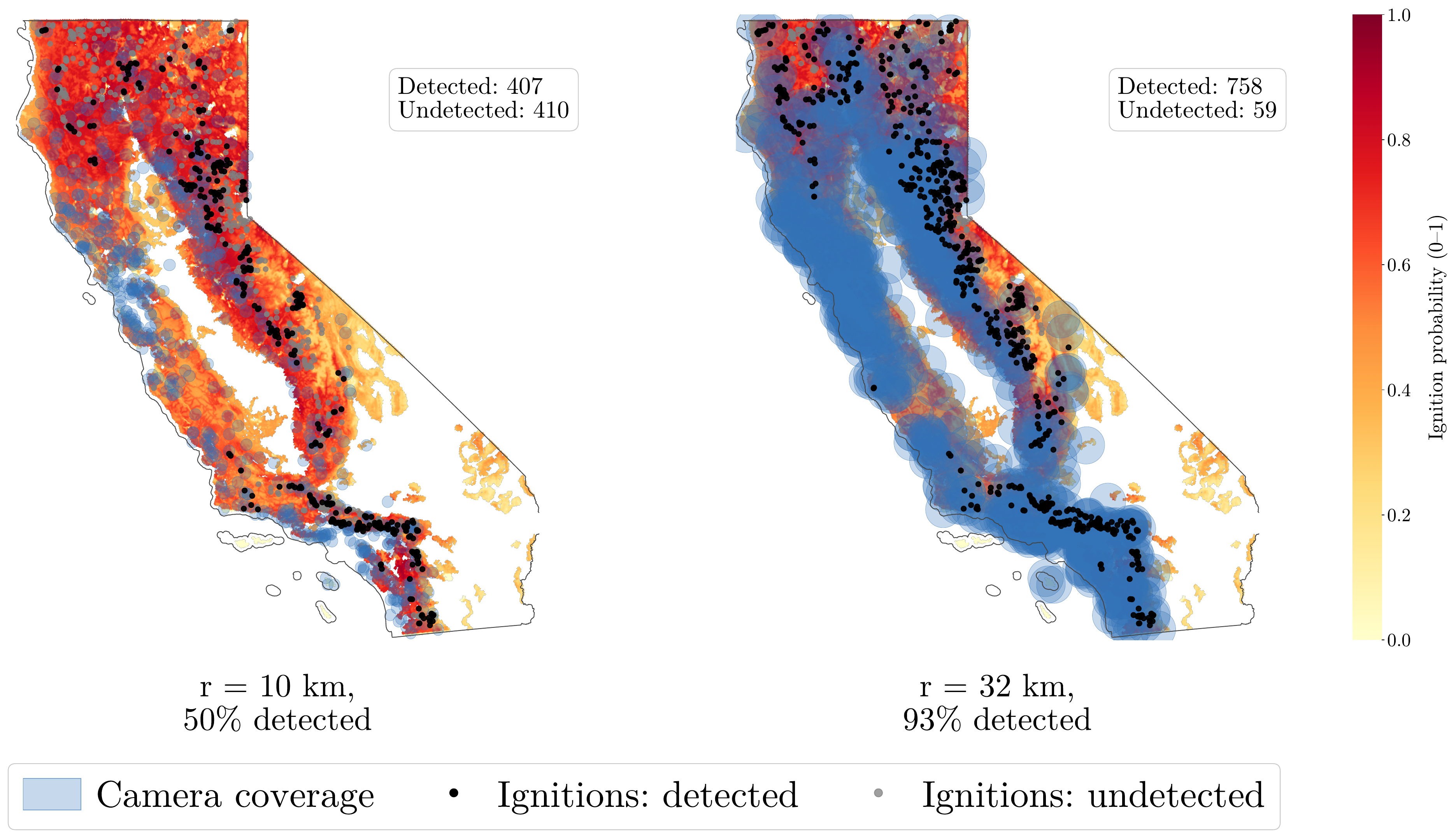}
    \caption{\textbf{ALERTCalifornia camera coverage of California 2024 wildfire ignitions ($n = 817$).}
Blue regions indicate the assumed camera detection range around each ALERTCalifornia site at 10~km and 32~km radii. Black dots denote detected ignitions and gray dots denote undetected ignitions. The 2024 dataset provides the cleanest out-of-sample evaluation, as these fires postdate the camera siting decisions.}
    \label{fig:alertcalifornia}
\end{figure*}

\danique{California's ALERTCalifornia program operates a statewide network of 699 fixed near-infrared camera sites on elevated terrain~\cite{alertcalifornia2026}. Rather than compare unpublished program costs~\cite{alertcalifornia2026}, we ask which drone budget matches camera coverage. On 2024 ignitions, the 10~km camera model (49.8\%) falls between the \$20 million (27.2\%) and \$50 million (57.9\%) drone networks.}
\danique{We use the camera network as a static-sensor baseline, evaluated on the same California ignition records (2021--2024; per-year $n$ in  Table~\ref{tab:alertcalifornia}). A fire is counted as detected if its ignition cell center lies within a Euclidean distance $r$ of any camera. This deliberately optimistic geometric model assumes perfect detection within radius, ignoring terrain occlusion, atmospheric attenuation, and the program's confirmation/ML-alerting layer. It therefore provides an upper bound on camera performance, and
any comparison in which the drone network exceeds it is correspondingly conservative. The effective detection range depends on terrain, weather, and fire characteristics, and ALERTCalifornia reports both a 60-mile line-of-sight range under ideal conditions and an online visualization radius of 20 miles. For our baseline experiments, we therefore
simulate detection rates across a range of assumed radii from 5 to 50~km (Table~\ref{tab:alertcalifornia}). Fig.~\ref{fig:alertcalifornia} illustrates coverage for $r=10$~km and $r=32$~km ($\sim$20 miles).}

\begin{table}[!h]
\caption{\textbf{ALERTCalifornia camera-network detection per year.} For each detection radius and year, we report the share of California wildfire ignitions whose location falls within the radius of at least one camera, with 95\% Wilson confidence intervals. Computed on the full California fire datasets (one column per year; 2024 is fully out-of-sample).}
\label{tab:alertcalifornia}
{\small
\begin{tabular}{@{}lcccc@{}}
\toprule
\textbf{Radius} & \textbf{2021} (n=981) & \textbf{2022} (n=903) & \textbf{2023} (n=992) & \textbf{2024} (n=817) \\
\midrule
5\,km & 19.6\% [17.2, 22.2] & 20.3\% [17.8, 23.0] & 19.2\% [16.8, 21.7] & 18.5\% [16.0, 21.3] \\
10\,km & 45.3\% [42.2, 48.4] & 52.4\% [49.1, 55.6] & 45.4\% [42.3, 48.5] & 49.8\% [46.4, 53.2] \\
20\,km & 78.8\% [76.1, 81.2] & 83.8\% [81.3, 86.1] & 75.9\% [73.1, 78.5] & 83.5\% [80.8, 85.9] \\
32\,km & 90.7\% [88.7, 92.4] & 95.1\% [93.5, 96.4] & 91.7\% [89.9, 93.3] & 92.8\% [90.8, 94.4] \\
50\,km & 100\% [99.6, 100] & 100\% [99.6, 100] & 100\% [99.6, 100] & 100\% [99.5, 100] \\
\bottomrule
\end{tabular}
}
\end{table}

\vspace{0.4cm}
\paragraph{Out-of-sample evaluation and ALERTCalifornia leakage.}
The ALERTCalifornia camera network has been deployed and progressively expanded
between 2017 and 2024, with camera siting partially informed by historical
ignition records that include fires from 2021--2023~\cite{alertcalifornia_history,alertcalifornia2026}.
This means that the 2021--2023 columns of Table~\ref{tab:alertcalifornia} probe
the network \emph{partly in-sample} from a siting-decision standpoint. The 2024
column provides the cleanest out-of-sample evaluation, given that the 2024
ignition data were not available when most camera placement decisions were made,
and we therefore use the 2024 dataset for Fig.~\ref{fig:alertcalifornia}.


On the fully out-of-sample 2024 ignitions ($n = 817$), the geometric camera model detects 49.8\% of fires at a 10~km radius and 92.8\% even at the optimistic 32~km sector range (Table~\ref{tab:alertcalifornia}). Because this model credits the cameras with perfect within-radius detection, these figures are upper bounds. On the same 2024 fires, the optimized drone network at \$100~million detects 97.2\% (95\% CI: 95.8--98.1\%), exceeding the most generous camera-coverage geometry. 



\section*{Discussion}

An optimized drone network operating at \$100 million detects 97.3\% of California's 2021--2024 wildfire ignitions, including 74\% within the first hour of ignition. Amortized over five years, that \$100 million budget is about \$20 million per year, or roughly 5\% of California's current annual expenditure on wildfire resource management and fire prevention \cite{lao_wildfire_faq_2025}. These results shift the policy question from whether rapid, near-complete wildfire detection is achievable to how monitoring investments should be structured to reach it.

The results reveal that spatial coverage is the dominant driver of detection performance. Once a fire falls within the reachable range of a charging station, detection rates exceed 95\% for the best-performing routing strategies across budget levels. Differences among these strategies primarily affect detection time rather than detection probability, with the best-performing strategy (TOP) achieving detection within the first hour for 74\% of fires at \$100 million. This indicates a transition from a coverage-limited regime, in which fires remain undetectable due to lack of reach, to a routing-limited regime, in which fires are reliably detected but with varying delays. Consequently, expanding the geographic footprint of the network emerges as the primary lever for improving detection, while investments in increasingly sophisticated routing yield diminishing returns.

Detection is not proportional to budget; it is threshold-driven. The relationship between budget and detection exhibits a pronounced nonlinearity, further underscoring the importance of coverage.  This behavior reflects the discrete nature of spatial coverage: a fire can only be detected if it lies within the reachable area of the network, while fires outside this footprint remain completely undetectable regardless of routing. Each additional charging station expands the reachable area by a fixed footprint, bringing previously unreachable fires into coverage. At \$20 million, 69\% of historical ignitions lie beyond the network’s reach; at \$100 million, this fraction decreases to 3\%. As a result, incremental investments below the coverage threshold yield limited value, since most fires remain fundamentally unreachable. These findings indicate that effective wildfire monitoring requires sufficient upfront investment to achieve near-complete coverage, rather than gradual scaling.

Under current hardware costs, drone-based monitoring is substantially more cost-effective per unit of covered area than static ground sensors. As a result, the placement optimization favors drone deployments across budget levels below \$500 million, with sensors not selected under these cost assumptions. However, the sensitivity analysis reveals that this ranking is not universal: as the relative cost of ground sensors decreases, a threshold emerges beyond which mixed deployments of drones and sensors become optimal, with sensors used to complement drone coverage in selected regions. While this transition does not occur within the budget and cost ranges examined in the baseline scenarios, it highlights that the relative roles of drones and sensors are sensitive to technology costs. Consequently, decisions regarding sensor network expansion should be evaluated against drone-based alternatives under realistic cost assumptions.

Our study has some limitations to be noted. We model detection as deterministic within what we call the certain-detection radius: once a fire is within this radius of a drone, detection is treated as both certain and instantaneous. In reality, detection probability decreases with distance from the fire. The true sensing footprint therefore has a reliable core, where the probability is close to 1, surrounded by a band in which detection is possible but not assured. Rather than model this gradient, we define the certain-detection radius as that inner core and discard the outer band. The resulting radius sits well below the nominal sensing ranges reported for comparable hardware (see Methods). This is what makes the assumption conservative. We count as guaranteed only those detections that would reliably occur, and ignore genuine coverage in the uncertain band. The opposite choice, treating the full nominal radius as certain, would credit detections in that outer band that would not reliably happen, and would therefore over-count. The detection rate, our headline result, is robust to this choice: we find that the bottleneck for detection is coverage, that is, whether a drone can reach a fire at all, and not per-pass detection success.


Our simulation also advances in discrete time steps at hourly granularity, so detection delays are binned hourly and we do not model sub-hour detection times. Hardware parameters including battery capacity, flight speed, and communication range are treated as fixed. The \$100 million budget is five-year hardware cost for stations and aircraft (Methods), not the cost of a permitted, staffed program. At that budget the network has 200 charging stations and 1{,}400 drones. Continuous beyond-visual-line-of-sight patrol at that scale is not authorized in the United States today. What does exist are smaller, legal deployments: California utilities already fly drone-in-a-box inspections under FAA Part 107 BVLOS waivers~\cite{pge_faa_bvlos_2026}, and commercial docking systems already run unattended battery swaps~\cite{Microavia2026}. Those programs show the hardware concept is in use, but they operate tens of docks on known assets, not a statewide search fleet. Staffing, waivers, and coordination with manned fire traffic would add a substantial multiple to the hardware figure we report. We leave that program cost out of the model. Land access, permitting for charging-station sites, and grid interconnection are likewise outside the model. 

\danique{Additionally, although our motivation is partly driven by climate-linked wildfire intensification, the Pyrologix risk map used in this study is static, calibrated to 2006–2020 ignitions with no explicit climate trend. Its continued predictive power through 2024 (Appendix A) is reassuring, but future climate-driven shifts in ignition patterns may require updated risk maps. This is not a limitation of the framework itself, which accepts any risk surface as input.}

\danique{The ALERTCalifornia comparison also has limitations. We assume a uniform effective detection radius for each camera, although performance varies with terrain and weather, and camera locations were selected using data that partly overlap the evaluation period. Accordingly, we present the camera comparison as a coverage-geometry
bound rather than an evaluation of the ALERTCalifornia program itself. A fully controlled comparison would require reconstructing camera detections from ALERTCalifornia's private confirmed-detection records. On the other hand, the drone network was optimized exclusively using pre-2021 inputs and hence evaluated entirely out-of-sample, whereas the current ALERTCalifornia camera network has evolved over time using historical fire records that include the evaluation period. Our results suggest that in regions with existing camera networks, optimization-based placement could guide future expansion, whereas in regions without such infrastructure, drone-based monitoring offers a more cost-effective foundation.}
The framework naturally supports additional policy objectives beyond maximizing wildfire detection. For example, fairness constraints could be incorporated to ensure equitable monitoring across communities while maintaining high detection performance. More broadly, the framework requires a spatial wildfire risk map, historical ignition locations, and a burnable-land mask as inputs, all of which are publicly available for most fire-prone regions globally. The results reveal structural insights that are expected to generalize across geographies: detection performance is primarily governed by spatial coverage, exhibiting a threshold behavior in budget; routing plays a secondary role once coverage is achieved; and the relative effectiveness of drones and static sensors depends on their cost, with drones dominating under current assumptions. While specific thresholds depend on regional topography, vegetation, and infrastructure costs, the overall patterns are robust, making applications to regions such as Mediterranean Europe, southeastern Australia, and the Brazilian Cerrado natural extensions. To facilitate further exploration, the optimization framework and datasets are released as open-source tools.

Several implications for decision-makers follow. Monitoring investments should target budget levels sufficient to achieve near-complete spatial coverage, rather than scaling incrementally below this threshold. Proposed expansions of static sensor networks should be evaluated against drone-based alternatives, given the substantial differences in cost-effectiveness. Regions with existing camera networks can benefit from optimization-based placement to guide expansion, whereas regions without such infrastructure may prioritize drone networks as the primary monitoring layer.

\section*{Methods}

\subsection*{Data sources and preprocessing}

\begin{figure*}[!h]
    \centering

    \begin{subfigure}{0.48\textwidth}
        \centering
        \includegraphics[width=\linewidth,height=0.31\textheight,keepaspectratio]{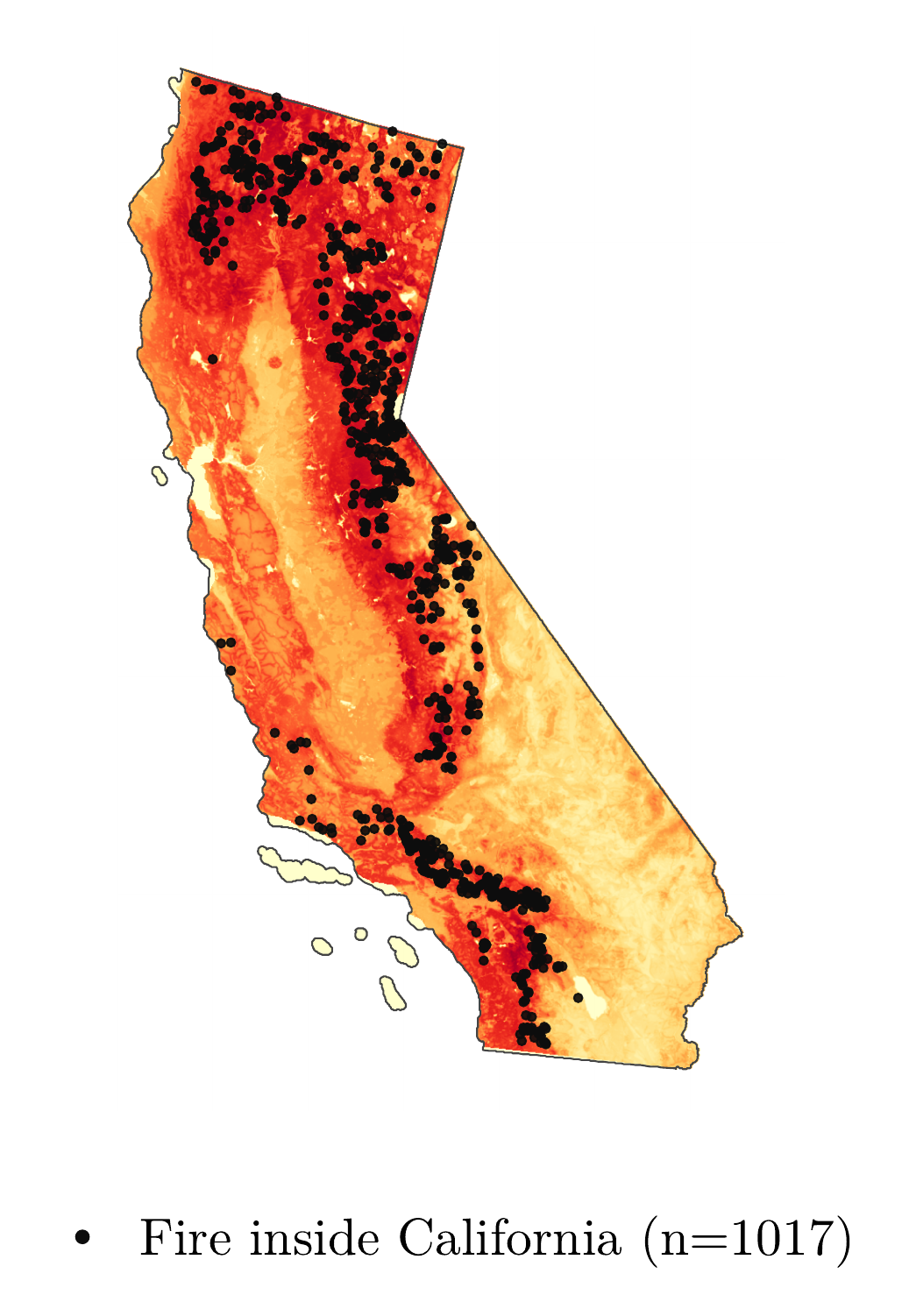}
        \caption{Study region}
        \label{fig:img1}
    \end{subfigure}
    \hfill
    \begin{subfigure}{0.48\textwidth}
        \centering
        \includegraphics[width=\linewidth,height=0.31\textheight,keepaspectratio]{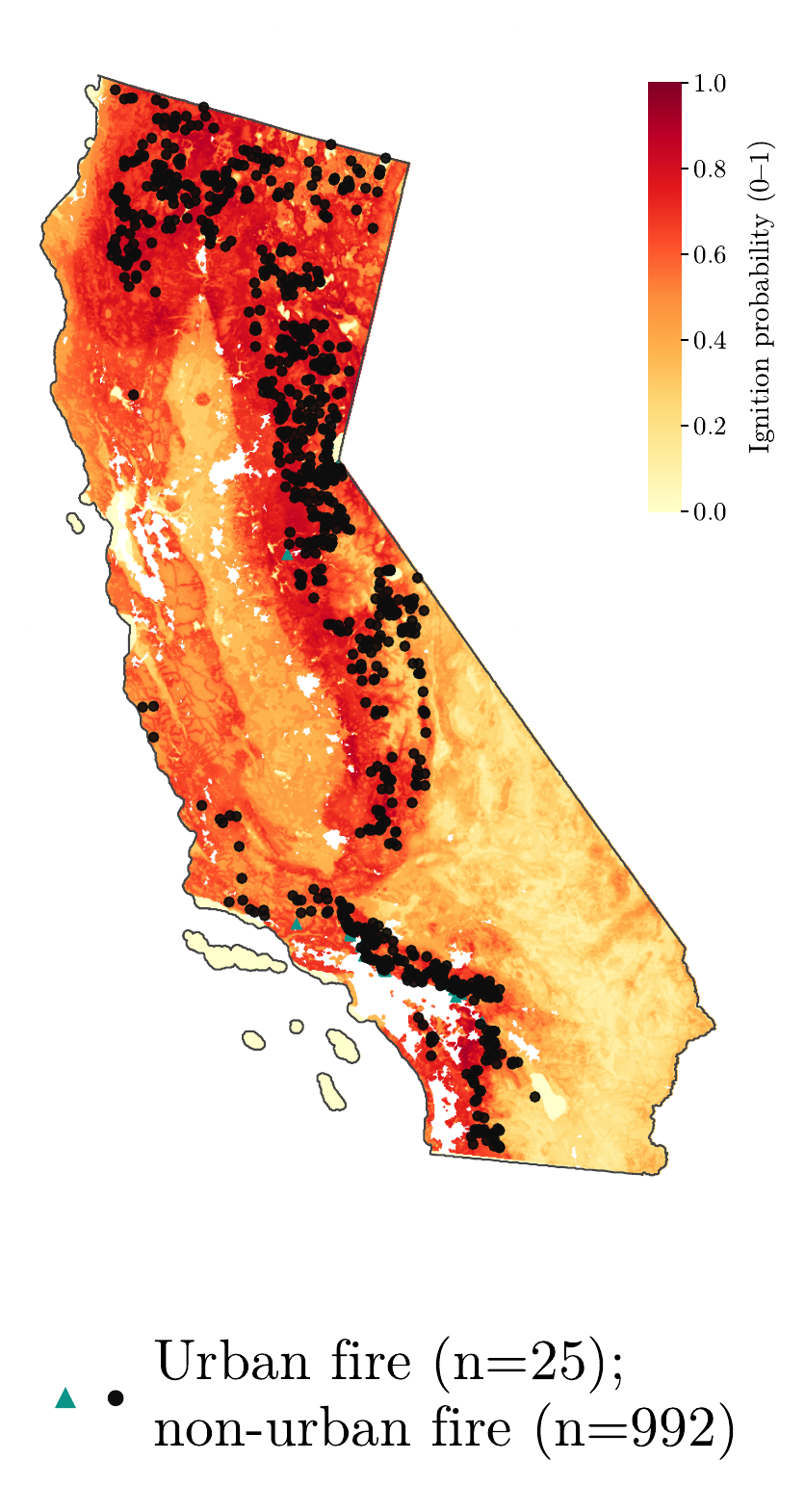}
        \caption{Excluding urban areas}
        \label{fig:img2}
    \end{subfigure}

    \par\medskip

    \begin{subfigure}{0.48\textwidth}
        \centering
        \includegraphics[width=\linewidth,height=0.31\textheight,keepaspectratio]{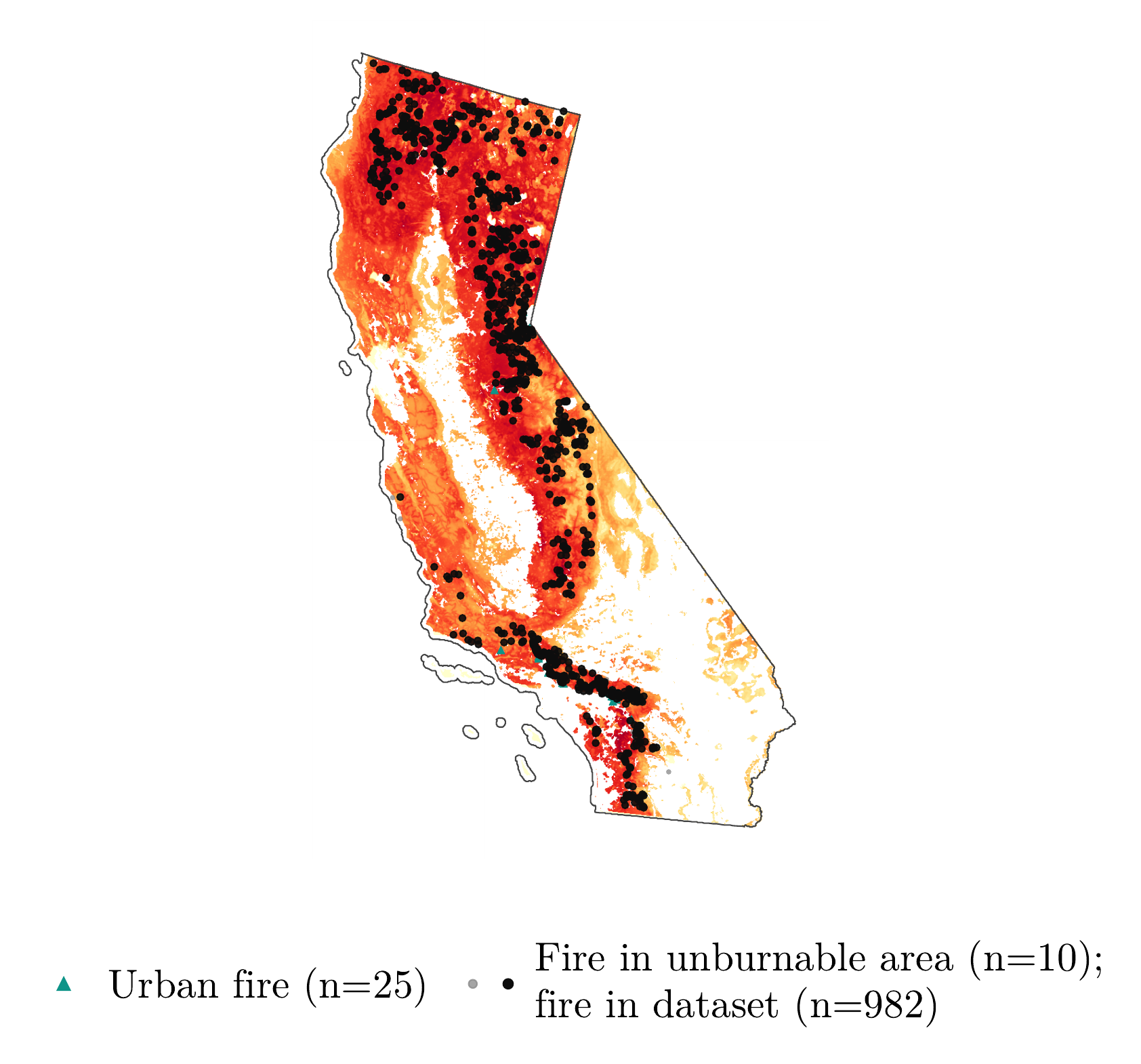}
        \caption{{\small Excluding unburnable areas}}
        \label{fig:img4}
    \end{subfigure}
    \hfill
    \begin{subfigure}{0.48\textwidth}
        \centering
        \includegraphics[width=\linewidth,height=0.31\textheight,keepaspectratio]{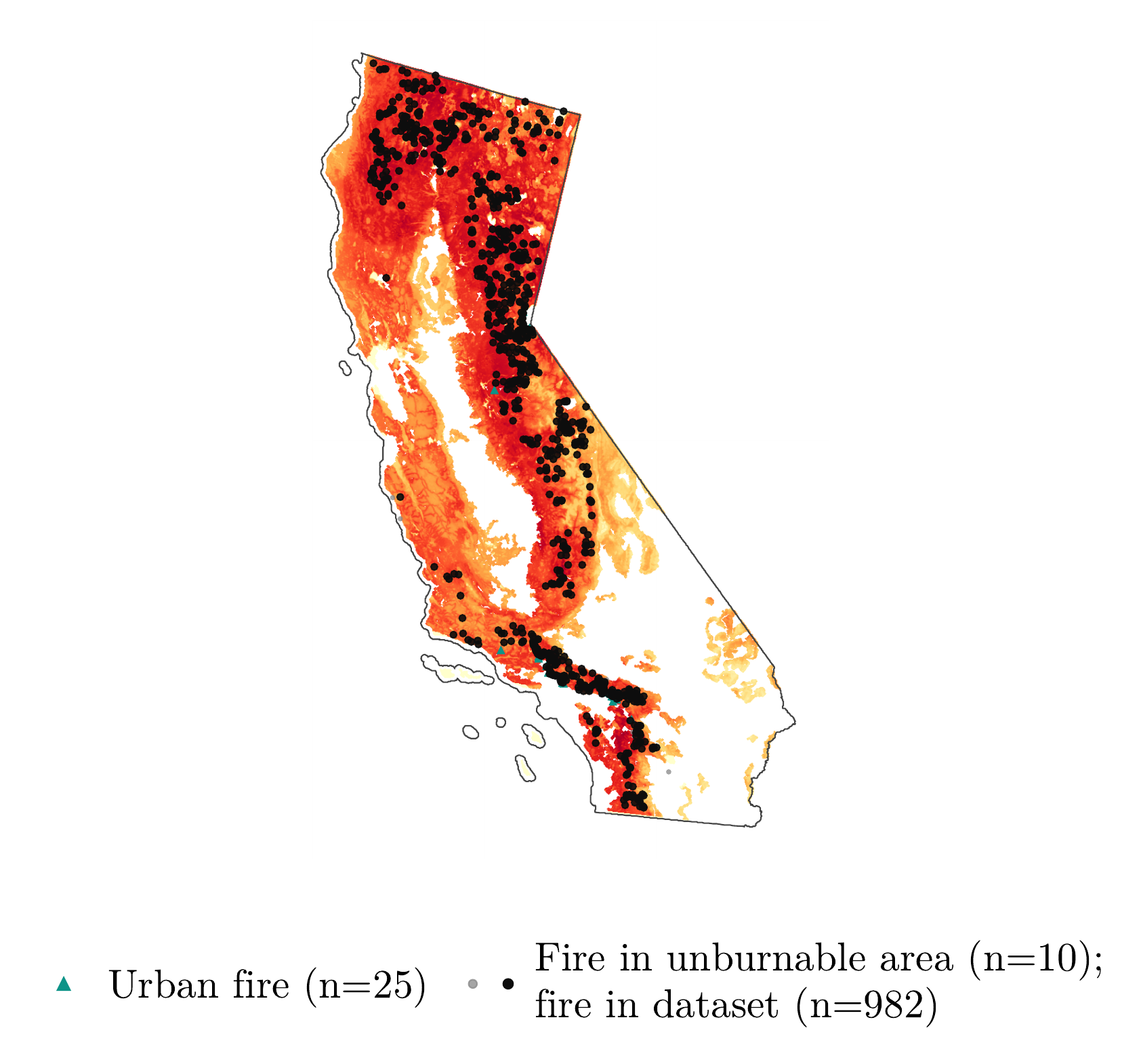}
        \caption{{\small Filtering small isolated regions}.}
        \label{fig:img5}
    \end{subfigure}
    \caption{\textbf{California case study for wildfire monitoring.}
Sequential preprocessing of the study region used for the wildfire monitoring experiments.
\textbf{(a)}~Geographic overview of California with the wildfire risk map and historical ignition points.
\textbf{(b)}~Burnable landscape after excluding urban areas, with removed ignition points highlighted.
\textbf{(c)}~Further refinement after excluding permanently unburnable regions, such as water bodies, deserts, and barren land, with additionally excluded ignition points highlighted.
\textbf{(d)}~Final study region after removing unconnected burnable regions smaller than $9 \times 9$ km.
    Illustrated for 2021 ignitions.}
    \label{fig:data}
\end{figure*}

To construct a realistic California case study, we integrate publicly available geospatial datasets, including wildfire risk estimates, historical ignition records, and land-use information, see  Fig.~\ref{fig:data}. 

The static wildfire risk map is obtained from Pyrologix and represents long-term wildfire ignition risk with probability values calibrated to observed annual ignition rates from 2006–2020~\cite{pyrologix}. Historical ignition points spanning 2021--2024 are obtained from the USFS fire occurrence dataset~\cite{usfs_ignition_points}. Because the risk map is calibrated exclusively to 2006--2020 ignition data, all benchmark years (2021--2024) are fully held out from the risk surface used for placement and routing. The evaluation is therefore strictly out-of-sample in time. Notably, a comparison with the California state boundary as defined by the U.S. Census map reveals that 30 ignition points across the four years (13 in 2021, 6 in 2022, 2 in 2023, and 9 in 2024) fall outside the official California boundary, highlighting minor spatial inconsistencies across public geospatial data sources. These ignition points are therefore excluded from the analysis.

We define a burnable landscape for California by excluding urban areas obtained from the U.S. Census Urban Areas 2020 dataset \cite{USCensusUrban2020} and permanently unburnable land-cover classes, including water bodies and barren land, identified through the Wildland Fire Potential Index (WFPI) 2020 dataset \cite{WFPI2020}. Only areas that are consistently unburnable throughout the year are excluded, thereby retaining regions that may be temporarily non-burnable because of seasonal effects such as snow cover. Unconnected burnable regions smaller than $9 \times 9$~km are then removed as noise.

To preserve ignition points near boundaries, a one-pixel buffer (1~km) is added around excluded areas. To avoid temporal leakage, we use the WFPI 2020 dataset rather than more recent releases.


\subsection*{Spatiotemporal grid construction}

The wildfire risk and historical ignition points data are provided on a projected equal-area grid at a spatial resolution of $d$ meters per cell and at an hourly temporal resolution. To make infrastructure placement and drone routing computationally tractable, we introduce a coarser \textit{operational grid} that aligns spatial and temporal scales with the physical capabilities of the monitoring hardware.

Each drone carries a sensor with detection radius $r$ meters, modeled as a square footprint of side $2r$ meters. We define the operational cell width as the nearest odd integer $c$ to $2r/d$ data cells, so that one operational cell corresponds to exactly one drone footprint. The data grid of $N_\mathrm{data} \times M_\mathrm{data}$ cells is tiled into non-overlapping $c \times c$ blocks, yielding an operational grid of
\[
N = \bigl\lfloor N_\mathrm{data}/c \bigr\rfloor \quad \text{and} \quad M = \bigl\lfloor M_\mathrm{data}/c \bigr\rfloor
\]
cells, on which all placement and routing decisions are made. The placement and routing models only act on this operational grid and do not consider the original data grid scale. A drone positioned at operational cell $(i,j)$ is assumed to observe the entire $c \times c$ block of data cells centered at that location. Planning at finer resolution would convey little additional operational information while scaling the decision problem quadratically.

The one-hour data timestep is both a constraint imposed by the temporal resolution of the historical ignition points data and a natural operational unit, since the drone battery capacity is set to one hour of flight time. A drone traveling at speed $v$ meters per minute can traverse
\[
k = \bigl\lfloor 60\,v / (c \cdot d) \bigr\rfloor
\]
operational cells per hour. Each such traversal constitutes one \textit{operational substep}, so a one-hour battery sustains exactly $k$ routing actions before requiring recharging. The simulation unfolds over hourly data timesteps, each subdivided into $k$ substeps during which drones move and fire detection is checked.

The wildfire risk map, originally at data resolution, is coarsened to the operational grid by averaging risk values within each $c \times c$ block. Each hourly risk value is then divided by $k$ and replicated $k$ times along the time axis, yielding per-substep risk values under the assumption of uniform ignition probability within each hour. All numerical parameter values used in the California case study are reported in Table~\ref{tab:benchmark_params}.


\subsection*{Infrastructure placement and drone allocation optimization}
The notation used throughout this formulation is summarized in Table~\ref{tab:infra_overview}.

\begin{table}[htbp]
\centering
\small
\caption{Overview of sets, parameters, and decision variables used in the infrastructure placement optimization model.}
\label{tab:infra_overview}
\begin{tabular}{lll}
\toprule
\textbf{Sets} & \textbf{Description} & \textbf{Definition} \\
\midrule
$\set I$ & Set of all grid points in the study region & $\{1,\dots,N\}\times\{1,\dots,M\}$ \\
$\set I_{g}$ & Feasible grid points for ground sensors & $\subseteq \set I$ \\
$\set I_{c}$ & Feasible grid points for charging stations & $\subseteq \set I$ \\
$\set P$ & Feasible (cell, station) pairs & $\subseteq \set I\times\set I_c$ \\
$\mathcal T$ & Set of time periods in the optimization horizon & $\{1,\dots,T\}$ \\
$\mathcal S$ & Set of drones & $\{1,\dots,S\}$ \\
$\mathcal C$ & Set of charging stations & $\{1,\dots,C\}$ \\
$\mathcal I_d$ & Set of feasible grid cells for drone movement & $\subseteq \mathcal I$ \\
$\mathcal V$ & Routing graph nodes (TOP) & $\mathcal I_{\mathrm{d}} \cup \{o,e\}$ \\
\midrule
\textbf{Parameters} & \multicolumn{2}{l}{\textbf{Description}} \\
\midrule
$r_{ij}$ & \multicolumn{2}{l}{Static wildfire risk score at grid point $(i,j)$} \\
$c_g$ & \multicolumn{2}{l}{Cost of one ground sensing station} \\
$c_c$ & \multicolumn{2}{l}{Cost of one charging station} \\
$c_d$ & \multicolumn{2}{l}{Cost of one drone} \\
$B$ & \multicolumn{2}{l}{Total deployment budget} \\
$\bar n_D$ & \multicolumn{2}{l}{Maximum number of drones assignable to one charging station} \\
$\epsilon$ & \multicolumn{2}{l}{Small regularization weight on infrastructure cost in the objective} \\
$\Delta_{s,k}$ & \multicolumn{2}{l}{Greedy incremental coverage fraction from the $k$-th drone at station $s$} \\
$r_{t+\delta,i}$ & \multicolumn{2}{l}{Wildfire risk at cell $i$ at global time $t+\delta$} \\
$B_{\max}$ & \multicolumn{2}{l}{Maximum drone battery capacity} \\
$\delta$ & \multicolumn{2}{l}{Time offset} \\
$T_{\mathrm{reev}}$ & \multicolumn{2}{l}{Number of periods executed before re-optimization} \\

$c_{ij}$ & \multicolumn{2}{l}{Travel cost between nodes $i$ and $j$ (TOP)} \\
$o$ & Origin depot (TOP) \\ 
$e$ & End depot (TOP) \\
\midrule
\textbf{Variables} & \multicolumn{2}{l}{\textbf{Description}} \\
\midrule
$x^g_{ij} \in \{0,1\}$ & \multicolumn{2}{l}{1 if a ground sensing station is placed at grid point $(i,j)$} \\
$x^c_{ij} \in \{0,1\}$ & \multicolumn{2}{l}{1 if a charging station is placed at grid point $(i,j)$} \\
$z_{ij,k} \in \{0,1\}$ & \multicolumn{2}{l}{1 if charging station $(i,j)$ has at least $k$ drones, $k=1,\dots,\bar n_D$} \\
$u_{ij,s} \in \{0,1\}$ & \multicolumn{2}{l}{1 if cell $(i,j)$ is assigned to charging station $s$ for coverage} \\
$w_{ij,s} \in [0,1]$ & \multicolumn{2}{l}{Coverage fraction contributed by station $s$ to cell $(i,j)$} \\
$\theta_{ij} \in [0,1]$ & \multicolumn{2}{l}{Infrastructure coverage level at grid point $(i,j)$} \\
$a_{its}\in\{0,1\}$ & \multicolumn{2}{l}{Drone $s$ at grid cell $i$ at time $t$} \\
$c_{jts}\in\{0,1\}$ & \multicolumn{2}{l}{Drone $s$ charging at station $j$ at time $t$} \\
$b_{ts}\in\mathbb Z_+$ & \multicolumn{2}{l}{Battery level} \\
$\zeta_{it}\in\mathbb Z_+$ & \multicolumn{2}{l}{Visit-history variable (min detection time)} \\
$w_{itt'}\in\{0,1\}$ & \multicolumn{2}{l}{Undetected risk indicator} \\
$y_{is}\in\{0,1\}$ & \multicolumn{2}{l}{Drone $s$ visits node $i$ (TOP)} \\
$x_{ijs}\in\{0,1\}$ & \multicolumn{2}{l}{Drone $s$ traverses arc $(i,j)$ (TOP)} \\
\bottomrule
\end{tabular}
\normalsize
\end{table}

\vspace{0.3cm}
\noindent \textit{Objective.} Given a budget, we jointly optimize the placement of ground sensors and charging stations, along with drone allocation, to maximize coverage of high-risk wildfire regions. Our objective is given by
\begin{align}
    \max_{\bm x^g,\bm x^c,\bm z,\bm u,\bm w,\bm \theta}\quad
& \sum_{(i,j)\in\mathcal I} r_{ij}\theta_{ij}
-\epsilon\left(
c_g\sum_{(i,j)\in\mathcal I_g} x^g_{ij}
+c_c\sum_{(i,j)\in\mathcal I_{\mathrm{c}}}x^c_{ij}
+c_d\sum_{(i,j)\in\mathcal I_{\mathrm{c}}}\sum_{k=1}^{\bar n_D}z_{ij,k}
\right).
\label{eq:infra_obj}
\end{align}
\noindent Here, we add a regularization term that ensures we avoid arbitrary over-deployment when the budget is nonbinding or weakly binding.\\
\\
\noindent \textit{Budget and resource allocation constraints.} These constraints are given by
\begin{align}
& c_g\sum_{(i,j)\in\mathcal I_g} x^g_{ij}
+c_c\sum_{(i,j)\in\mathcal I_c}x^c_{ij}
+c_d\sum_{(i,j)\in\mathcal I_c}\sum_{k=1}^{\bar n_D}z_{ij,k}
\le B,
\label{eq:infra_budget}\\
& z_{ij,k}\le x^c_{ij},
&&  (i,j)\in\mathcal I_c,\; k=1,\dots,\bar n_D,
\label{eq:infra_link_station}\\
& z_{ij,k}\ge z_{ij,k+1},
&&  (i,j)\in\mathcal I_c,\; k=1,\dots,\bar n_D-1.
\label{eq:infra_link_order}
\end{align}
Constraint \eqref{eq:infra_budget} enforces the overall deployment budget. 
Constraint \eqref{eq:infra_link_station} ensures drones can only be allocated where a charging station is installed, and constraint \eqref{eq:infra_link_order} enforces the ordering so that $\sum_k z_{ij,k}$ correctly counts the number of drones at station $(i,j)$.\\
\\
\noindent \textit{Infrastructure siting constraint.} To prevent colocating both infrastructure types at the same feasible location, we impose
\begin{align}
& x^g_{ij}+x^c_{ij}\le 1,
&&  (i,j)\in\mathcal I_g \cap \mathcal I_c.
\label{eq:infra_exclusion}
\end{align}
Feasibility here is defined purely by the burnable-land mask (Fig. \ref{fig:data}). Restricting $\mathcal I_g$ and $\mathcal I_c$ with land-ownership, suitable terrain, and power access is a natural extension.\\
\\
\noindent \textit{Cell assignment and coverage constraints.}
The feasible set $\mathcal P \subseteq \mathcal I \times \mathcal I_c$ contains all (cell, station) pairs where the cell is reachable by a drone launched from the station within one battery charge, determined by breadth-first search on the masked operational grid. Each cell is assigned to at most one charging station:
\begin{align}
& \sum_{s:\,(i,j,s)\in\mathcal P} u_{ij,s}\le 1,
&&  (i,j)\in\mathcal I.
\label{eq:infra_assign}
\end{align}
The coverage fraction $w_{ij,s}$ that station $s$ contributes to cell $(i,j)$ is a piecewise-linear, diminishing-returns function of the drone count, precomputed with a greedy set-cover heuristic over drone patrol paths within the station's zone. For each additional drone level $k$, the heuristic selects the patrol path that maximizes newly covered risk; the resulting incremental coverage fraction $\Delta_{s,k}$ is then applied uniformly to all reachable cells in the station's zone. We set $\bar n_D=k$ (here 7), the drone count at which a station's reachable zone is fully covered; additional drones add cost but no coverage. This gives
\begin{align}
& w_{ij,s}\le \sum_{k=1}^{\bar n_D}\Delta_{s,k}\,z_{s,k},
&&  (i,j,s)\in\mathcal P,
\label{eq:infra_wbound_drone}\\
& w_{ij,s}\le u_{ij,s},
&&  (i,j,s)\in\mathcal P,
\label{eq:infra_wbound_assign}\\
& w_{ij,s}\ge 0,
&&  (i,j,s)\in\mathcal P.
\label{eq:infra_wpos}
\end{align}
Ground sensing stations provide direct local coverage at their installed grid cells. The overall coverage level at each cell combines ground sensing with drone-based surveillance:
\begin{align}
& \theta_{ij}\ge x^g_{ij},
&&  (i,j)\in\mathcal I_g,
\label{eq:infra_ground_only}\\
& \theta_{ij}\le x^g_{ij}+\sum_{s:\,(i,j,s)\in\mathcal P}w_{ij,s},
&&  (i,j)\in\mathcal I_g,
\label{eq:infra_coverage_prime}\\
& \theta_{ij}\le \sum_{s:\,(i,j,s)\in\mathcal P}w_{ij,s},
&&  (i,j)\in\mathcal I\setminus\mathcal I_g.
\label{eq:infra_coverage_nonprime}
\end{align}
Constraints \eqref{eq:infra_coverage_prime}--\eqref{eq:infra_coverage_nonprime} combine direct sensing from ground stations with indirect drone-based coverage from charging stations.\\
\\
\noindent \textit{Variable domain constraints.} Finally, the infrastructure decisions satisfy
\begin{align}
& 0\le \theta_{ij}\le 1,
&&  (i,j)\in\mathcal I,
\label{eq:infra_theta}\\
& x^g_{ij}\in\{0,1\},
&&  (i,j)\in\mathcal I_g,
\label{eq:infra_binary_g}\\
& x^c_{ij}\in\{0,1\},
&&  (i,j)\in\mathcal I_c,
\label{eq:infra_binary_c}\\
& z_{ij,k}\in\{0,1\},
&&  (i,j)\in\mathcal I_c,\; k=1,\dots,\bar n_D,
\label{eq:infra_binary_z}\\
& u_{ij,s}\in\{0,1\},
&&  (i,j,s)\in\mathcal P.
\label{eq:infra_binary_u}
\end{align}


\subsection*{Drone routing optimization}
We propose three complementary drone-routing strategies for proactive wildfire monitoring: a Team Orienteering Problem-based strategy (TOP),  a maximum-coverage strategy (MaxCov), and a minimum-detection-time strategy (LinearMinTime). TOP builds on the classical Team Orienteering Problem introduced by Chao et al.~\cite{chao1996team}. MaxCov is adopted from~\cite{puech2026wfdronebench} and builds on the classical maximal covering location problem introduced by Church and ReVelle~\cite{MCLPChurchRevelle}. LinearMinTime is introduced here to directly prioritize reduction of wildfire detection delay. While the strategies share a common routing and feasibility structure, they differ in how surveillance performance is quantified and optimized.

The sets, parameters, and decision variables of the routing models are summarized in Table~\ref{tab:infra_overview}. To restrict operations to accessible regions, we apply a binary spatial mask over the grid, excluding infeasible areas such as water bodies or terrain obstacles. A breadth-first search (BFS) from charging-station locations is then used to identify all reachable cells, defining the feasible set $\set I_d$. This ensures that routing is confined to connected, operationally feasible regions while reducing computational complexity.

\subsubsection*{Team Orienteering Problem}

\noindent \textit{Objective.}
The TOP strategy formulates wildfire monitoring as a route-selection problem in which drones seek to maximize the cumulative wildfire risk observed within a finite battery horizon. Each drone departs from a charging station, visits a subset of feasible grid cells, and must return to a charging station before battery depletion.  The objective is given by
\begin{equation}
\label{eq:top_obj}
\max \quad
\sum_{s \in \mathcal{S}}
\sum_{i \in \mathcal{I}_{\mathrm{d}}}
r_i y_{is}.
\end{equation}
Unlike time-indexed routing formulations, TOP optimizes complete drone routes directly on a routing graph rather than sequential drone positions over time.\\
\\
TOP plans a complete route for each drone over a single battery horizon, and the problem is re-solved once a drone completes its battery cycle. The risk values $r_i$ are updated between these successive solves to reflect the mitigating effect of active monitoring. When a drone surveils a cell, that cell's accumulated ignition risk is reset to zero; while a cell remains unmonitored, its risk instead grows back linearly, accruing the cell's baseline per-step ignition probability $r_0$ at each step, so that after $k$ unmonitored steps its risk equals $k\,r_0$. This linear accumulation is the first-order Taylor approximation of the exact probability that at least one ignition has occurred over $k$ independent steps, $1-(1-r_0)^k = k\,r_0 + O\!\big((k r_0)^2\big)$, which is accurate here because per-cell hourly ignition probabilities are small.\\
\\
\noindent \textit{Artificial origin and destination depots.}
The artificial depot nodes, $o$ and $e$, correspond to the same physical charging station, but  are modeled as distinct nodes to separate departure and return flows, a standard approach in vehicle-routing formulations.\\
\\
\noindent \textit{Routing constraints.}
The constraints are given by
\begin{align}
& \sum_{s \in \mathcal{S}} y_{is} \le 1,
&& i \in \mathcal{I}_{\mathrm{d}},
\label{eq:top_unique_visit}
\\
& \sum_{j \in \mathcal{V}^{-}} x_{ojs} = 1,
&& s \in \mathcal{S},
\label{eq:top_start}
\\
& \sum_{i \in \mathcal{V}^{+}} x_{ies} = 1,
&& s \in \mathcal{S},
\label{eq:top_end}
\\
& x_{jos}=0,
&& j \in \mathcal{V},\ s \in \mathcal{S},
\label{eq:top_no_enter_origin}
\\
& x_{ejs}=0,
&& j \in \mathcal{V},\ s \in \mathcal{S},
\label{eq:top_no_leave_dest}
\\
& \sum_{j \in \mathcal{V}^{-}\setminus\{i\}} x_{ijs}
=
y_{is},
&& i \in \mathcal{I}_{\mathrm{d}},\ s \in \mathcal{S},
\label{eq:top_outflow}
\\
& \sum_{j \in \mathcal{V}^{+}\setminus\{i\}} x_{jis}
=
y_{is},
&& i \in \mathcal{I}_{\mathrm{d}},\ s \in \mathcal{S},
\label{eq:top_inflow}
\end{align}
where $\mathcal{V}=\mathcal{I}_{\mathrm{d}}\cup\{o,e\}$, $\mathcal{V}^{+}=\mathcal{I}_{\mathrm{d}}\cup\{o\}$, and $\mathcal{V}^{-}=\mathcal{I}_{\mathrm{d}}\cup\{e\}$.

Constraints \eqref{eq:top_unique_visit} ensure that each grid point is served by at most one drone, thereby avoiding duplicate reward collection across the fleet. Constraints \eqref{eq:top_start} and \eqref{eq:top_end} ensure that each drone departs from the charging station and returns exactly once. Constraints \eqref{eq:top_outflow} and \eqref{eq:top_inflow} ensure that if a cell is visited, it must be entered and exited exactly once.\\
\\
\noindent \textit{Battery and travel constraints.}
The battery feasibility constraints are given by
\begin{align}
& \sum_{i \in \mathcal{V}^{+}}
\sum_{\substack{j \in \mathcal{V}^{-}\\ j \neq i}}
c_{ij}x_{ijs}
\le B_{\text{max}},
&& s \in \mathcal{S},
\label{eq:top_battery}
\end{align}
\noindent ensuring that the total travel distance of each drone remains within the battery capacity $B_{\text{max}}$. The travel cost $c_{ij}$ is defined using the $L_{\infty}$ metric.
\\
\\
\textit{Symmetry-breaking constraints.} Constraints \eqref{eq:top_symmetry} improve computational efficiency by reducing equivalent solutions that differ only by permutation of drone labels.
\begin{align}
    & \sum_{i \in \mathcal{I}_{\mathrm{d}}} r_i y_{i,s+1}
\le
\sum_{i \in \mathcal{I}_{\mathrm{d}}} r_i y_{is},
&& s=1,\dots,|\mathcal{S}|-1.
\label{eq:top_symmetry}
\end{align}
\\
\\
\noindent \textit{Particle swarm optimization (PSO) routing heuristic.}
The TOP formulation yields a large combinatorial problem whose complexity scales rapidly with the number of drones, grid cells, charging stations, and arcs, making exact optimization intractable for large or time-sensitive wildfire applications. We therefore adopt the PSO heuristic in~\cite{PSO} to efficiently explore candidate multi-drone routes and obtain high-quality solutions at lower cost. The reported TOP-based routing results in Table~\ref{tab:detection} are thus heuristic. Path connectivity is enforced by construction in the PSO heuristic; the displayed model does not include subtour-elimination constraints. We adapt this algorithm to our grid structure to make it substantially faster, principally through a sparse split procedure that exploits the small number of charging stations, conservative filtering of local-search moves that cannot change the solution, and incremental route updates; 
Appendix B details these changes, and the implementation is available in the code released with this paper (see the Code Availability statement).

\subsubsection*{Maximizing coverage}
We adopt the risk-aware maximum-coverage routing strategy proposed in~\citep{puech2026wfdronebench}. The strategy uses a rolling-horizon mixed-integer optimization model that repeatedly routes drones to maximize the cumulative wildfire risk of newly monitored areas while accounting for operational constraints such as battery endurance, charging requirements, communication limits, and feasible drone movements. In contrast to TOP, which is simply re-solved once per battery cycle, MaxCov is a true rolling-horizon scheme: it is re-optimized every $T_{\text{reev}}$ steps within the planning horizon, allowing surveillance to adapt more responsively as conditions evolve. The wildfire risk map is updated at each re-optimization to reflect the mitigating effect of active monitoring, following the same mechanism as in the TOP strategy: a surveilled cell has its accumulated ignition risk reset to zero, while an unmonitored cell regrows its risk linearly at the baseline per-step rate $r_0$. Full methodological details are provided in~\citep{puech2026wfdronebench}.

\subsubsection*{Minimizing detection time} 
Similar to the maximum-coverage approach, this strategy is implemented within a rolling-horizon framework.\\
\\
\noindent \textit{Objective.}
The goal is to minimize the expected wildfire detection time, defined as the time between ignition and first detection by any monitoring device.

Detection time is stochastic, as it depends on both the random ignition process and operational decisions such as sensor placement and drone routing. An exact formulation leads to a nonlinear stochastic integer program, since detection probabilities depend on the full history of visits to each cell and the cumulative probability that ignition has occurred since the most recent visit.

To obtain a tractable model, we use a linear surrogate objective. Instead of minimizing expected detection time directly, we minimize cumulative undetected risk, defined as the probability that a cell is burning between successive visits. This proxy captures detection delay: risk accumulates the longer a high-risk cell remains unvisited, encouraging earlier and more frequent monitoring of such regions.



\begin{equation}
\label{eq:routing_obj_paper}
\min_{\bm \zeta, \bm w} \quad 
\sum_{t \in \mathcal{T}}
\sum_{t'=1}^{t}
\sum_{i \in \mathcal{I}_{\text{d}}}
\left(
w_{itt'}\, r_{t'+\delta,i}
\right),
\end{equation}
where the summation
\[
\sum_{t'=1}^{t} w_{itt'} r_{t'+\delta,i}
\]
captures the cumulative probability mass that a fire may remain undetected since the most recent visit. 
\\
\\
\textit{Constraints.} The drone movement, battery, and initialization constraints are the same as for the max-coverage strategy. In addition, we have the following visit-history and risk accumulation constraints
\begin{align}
& \zeta_{it} \ge \zeta_{i,t-1},
&& i \in \set I_{\text{d}},\ t \in \set T \setminus \{1\}
\label{eq:routing_zeta_monotone}
\\
& \zeta_{i1} \ge 0,
&& i \in \set I_{\text{d}},
\label{eq:routing_zeta_init_lb}
\\
& \zeta_{it} \le t,
&& i \in \set I_{\text{d}},\ t \in \mathcal{T},
\label{eq:routing_zeta_time}
\\
& \zeta_{it}
\le
\zeta_{i,t-1}
+
t\sum_{s\in\mathcal S} a_{its},
&& i \in \set I_{\text{d}},\ t \in \set T \setminus \{1\}
\label{eq:routing_zeta_visit}
\\
& \zeta_{i1}
\le
\sum_{s\in\mathcal S} a_{i1s},
&& i \in \set I_{\text{d}},
\label{eq:routing_zeta_first}
\\
& t' - \zeta_{it} \le t\, w_{itt'},
&& i \in \set I_{\text{d}},\ t \in \mathcal{T},\ t'=1,\dots,t,
\label{eq:routing_w_link}
\\
& w_{jtt'} = 1,
&& j \in \mathcal{C},\ t \in \mathcal{T},\ t'=1,\dots,t.
\label{eq:routing_w_charge}
\end{align}
Constraints \eqref{eq:routing_zeta_monotone}--\eqref{eq:routing_zeta_first} ensure that this visit-history variable is nondecreasing over time and is updated whenever a drone visits cell \(i\). 
Constraint \eqref{eq:routing_w_link} enforces that \(w_{itt'}=1\) unless cell \(i\) has been visited after time \(t'\). Constraint \eqref{eq:routing_w_charge} fixes charging station cells as permanently visited. Finally, the decision variables 
satisfy 
\begin{align}
& a_{its}\in\{0,1\},
&& i\in\mathcal I_d,\ t\in\mathcal T,\ s\in\mathcal S,
\label{eq:routing_binary_a}\\
& c_{jts}\in\{0,1\},
&& j\in\mathcal C,\ t\in\mathcal T,\ s\in\mathcal S,
\label{eq:routing_binary_c}\\
& b_{ts}\in\mathbb Z_+,
&& t\in\mathcal T,\ s\in\mathcal S.
\label{eq:routing_integer_b}\\
& w_{itt'} \in \{0,1\}, && i \in \set I_d, \ t, t' \in \set T, \\
& \zeta_{it} \in \mathbb{Z}_+, && i \in \set I_d, \ t \in \set T.
\label{eq:routing_w_dom_vec}
\end{align}

\subsection*{Cost and hardware parameterization}\label{sec:params}
Table~\ref{tab:benchmark_params} lists the numerical settings used in our experiments. \\

\begin{table}[t]
\caption{\textbf{Benchmark and drone parameters used in our experiments.}}
\label{tab:benchmark_params}
\begin{tabular}{@{}p{0.42\linewidth}p{0.48\linewidth}@{}}
\toprule
\textbf{Parameter} & \textbf{Value} \\
\midrule
\multicolumn{2}{@{}l}{\textit{Spatial-temporal grid}} \\
Data resolution ($d$) & $1\ \mathrm{km}$ per cell \\
Data grid size & $1309 \times 805$ cells \\
Data timestep & $1\ \mathrm{hour}$ \\
\midrule
\multicolumn{2}{@{}l}{\textit{Drone and motion model}} \\
Drone speed ($v$) & $600\ \mathrm{m\,min^{-1}}$ \\
Drone coverage radius ($r$) & $2900\ \mathrm{m}$ (operational cell $c = 5$ data cells) \\
Battery capacity & $1\ \mathrm{hour}$ ($k = 7$ operational substeps) \\
Transmission range & $50\ \mathrm{km}$ \\
Mixed-integer solver & Gurobi; time limit 2\,mins per routing step solve \\
Re-optimization interval ($T_{\mathrm{reev}}$) & $5$ operational substeps \\
\midrule
\multicolumn{2}{@{}l}{\textit{Budget optimization}} \\
Total budgets & USD $20$ / $50$ / $75$ / $100$ / $500$ million \\
Cost per ground sensor & USD $100,000$  \\
Cost per charging station & USD $150,000$ \\
Cost per drone & USD $50,000$ \\
Regularization weight ($\epsilon$) for the $\$500$ million run & $1$ \\
\bottomrule
\end{tabular}
\end{table}

\noindent \textit{Budget optimization parameters.}
Public wildfire-prevention funding spans a wide range, with programs allocating on the order of tens to hundreds of millions of dollars annually \cite{calfire_wildfire_grants_2025, CalMattersWildfireFunding2025}. Motivated by this scale, we consider budget levels of \$20M, \$50M, \$75M, \$100M, and \$500M.

\medskip
\noindent \textit{Infrastructure cost assumptions.}
The unit costs are assumed to be \$100{,}000 per ground sensor, \$150{,}000 per charging station, and \$50{,}000 per drone. These deployment-scale estimates represent the total cost of ownership over an approximately five-year horizon, combining initial capital expenditures with operating and maintenance (O\&M) costs. These unit costs are hardware total cost of ownership. They do not include airspace authorization or staffing costs. Because comprehensive, standardized pricing for these emerging systems is rarely published, we did not derive these values solely from the literature and calibrated them instead in consultation with a group of California-based researchers actively developing and deploying drones for wildfire detection. We cross-checked them against available public cost data, summarized below.

For ground sensors, we deliberately model a fully deployed, standalone fixed-camera station. Public figures indicate that the camera hardware itself is inexpensive, on the order of \$2{,}500--\$3{,}000, but that the supporting site infrastructure (tower or mast, solar power, backhaul communications, and installation) can add up to roughly \$75{,}000 depending on terrain and accessibility~\cite{GovTech2018WildfireCamera}, consistent with per-camera capital costs of \$15{,}000--\$100{,}000 reported for the Oregon network~\cite{KLCC2026OregonCameras}. Reported annual operating costs likewise span a wide range, from roughly \$6{,}000 per year where existing power and connectivity can be reused~\cite{GovTech2025SantaFeAI} to on the order of \$36{,}000 per station-year for large utility-scale deployments that include field maintenance and around-the-clock monitoring (e.g., Hawaiian Electric's \$14M, five-year contract for 78 stations)~\cite{MauiNews2024HawaiianElectric}. Combining installed capital cost with several years of O\&M places a representative standalone station near our \$100{,}000 figure, which sits at the conservative end of the reported range.

For charging stations, autonomous drone-in-a-box systems bundle the docking enclosure, robotic battery-management system, communications, and installation; the Microavia platform, for example, performs a fully robotic battery swap in roughly two minutes, enabling near-continuous operation with minimal human intervention~\cite{Microavia2026}. Public pricing for such industrial systems is largely project-specific and undisclosed, so we anchor the \$150{,}000 figure to the deployment experience gathered in this consultation, with vendor specifications corroborating the assumed operational capability. For drones, combining a long-range electric VTOL airframe~\cite{mugin_ev350} with a thermal payload, optical sensors, stabilization, and a professional autopilot stack~\cite{FLIRBoson640, SonyBlockCamera, GremsyIndustrialGimbal, CubePilotAutopilot} yields a mission-ready acquisition cost of roughly \$18{,}000--\$30{,}000; because sustained fire-season operation is dominated by recurring battery replacement and servicing rather than airframe capital, we adopt a higher effective five-year cost of \$50{,}000 per drone.

All three values are exposed as parameters in our open-source framework, and the cost-sensitivity analysis (Fig.~\ref{fig:breakeven_costsensitivity}) characterizes how the optimal infrastructure mix responds when they are varied.

\medskip
\noindent \textit{Operational drone parameters.}
The operational parameters are inspired by the Mugin EV350~\cite{mugin_ev350}, a long-range electric VTOL drone used for autonomous environmental monitoring, and, as with the cost assumptions above, were calibrated in consultation with the same group of researchers. We use a 1~km spatial grid covering $1309 \times 805$ cells with hourly time steps. 

Our model credits a fire as detected when it falls within a drone's \emph{certain-detection radius}. This is a conservative coverage abstraction rather than an assumption that sensors never miss fires: the true detection probability $P_d(d)$ decreases smoothly with the distance $d$ between fire and sensor and depends on fire area, sensor modality, and atmospheric conditions~\cite{WildfireDetectionDLKim,WildfireDetectionDLKumar24,WildfireDetectionDLYang}. We replace this surface by the indicator $\mathds{1}[d \le r^\ast]$ with $r^\ast = 2{,}900$~m, chosen well inside the regime where $P_d \approx 1$, i.e.\ substantially below the nominal sensing range of the reference hardware. Coupling the threshold to range in this way is more physically faithful than collapsing $P_d$ to a single range-independent scalar, since $P_d$ is intrinsically a function of range and fire size; we therefore absorb detection reliability into the coverage geometry and treat the boundary conservatively. Because we discard all genuine detections that occur at $d > r^\ast$, where $0 < P_d < 1$, this step-function approximation \emph{lower-bounds} the expected number of detections rather than inflating it.

Two mechanisms make the detection \emph{rate} insensitive to the residual per-pass miss probability at $d \le r^\ast$. First, fires persist and grow, so a small fire not detected on an early pass becomes easier to detect on subsequent passes, for which $r^\ast$ is increasingly conservative. Second, reachable cells are revisited repeatedly within the detection window, so the cumulative detection probability $1 - (1 - P_d)^{m}$ over $m$ passes approaches one even for moderate $P_d$. The binding constraint on the detection rate is therefore spatial reachability, a geometric property independent of $P_d$, consistent with our central finding that placement governs detection rate while routing governs detection speed.

We adopt a conservative one-hour effective battery life, reflecting reduced capacity under adverse conditions such as high winds or temperature extremes. We assume immediate battery swap at charging stations and no mid-air collisions between drones. The effective scanning speed of 600~m\,min$^{-1}$ trades nominal flight speed for the dwell time required for reliable detection, and was estimated from camera sensor characteristics and field experience by the same group of researchers. The 50~km transmission range is non-binding in practice: given the one-hour battery and the 600~m\,min$^{-1}$ cruising speed, a drone can cover at most 36~km round-trip from its charging station, well within communication range under all routing strategies evaluated.

All parameters are inputs to our open-source framework; we encourage readers to adjust them and re-run the benchmarks to reflect alternative drone platforms, cost structures, or detection assumptions.

\subsection*{Placement and routing computation times.}
Sensor placement at the \$75\,M budget was solved as a single mixed-integer program (Gurobi~12, 32 CPUs), terminating at the 1800 s wall-clock limit with a MIP gap of 0.01\,\%. Each charging station computes a 24-hour routing plan that is re-optimized every five control substeps and cached for reuse across all benchmark years. The reported runtimes were measured during this initial (cold-cache) computation and are normalized per simulated hour of routing. The mixed-integer routing formulations (MaxCov and LinearMinTime) use a 120\,s time limit for each Gurobi subproblem. Mean (median) solver times per simulated hour are 31.9 (28.4)\,s for TOP, 107.0 (106.7)\,s for MaxCov, and 184.2 (210.4)\,s for LinearMinTime. LinearMinTime exceeds 120\,s because a simulated hour comprises more than one subproblem. For MaxCov, 53\,\% of the per-step Gurobi sub-problems solve to proven optimality within the 120\,s time limit, while the remainder terminate close to optimality. TOP is substantially faster than the two mixed-integer formulations. Once the routing solutions have been computed, each additional year requires only fire-simulation replay with negligible routing overhead.

\subsection*{Statistical analysis}
Detection time is recorded at one-hour resolution. $\Delta t = k$ indicates that a fire was detected during the interval $[k,\, k+1)$ hours after ignition, yielding ordinal integer values $\Delta t \in \{0,1,2,3,4,5\}$. A fire not detected within six hours of ignition is counted as missed. The primary speed metric is the proportion of all fires with $\Delta t = 0$, i.e., detection occurring within the first hour. Confidence intervals for all proportions (detection rate, reachable rate, and proportion within first hour) use the 95\% Wilson score interval for binomial proportions. Differences in detection speed across routing strategies are assessed with the two-sided Mann-Whitney $U$ test applied to the ordinal $\Delta t$ values among jointly detected fires.

\section*{Acknowledgements}
The second author is funded by the Netherlands Organisation for Scientific Research (NWO) under the Rubicon grant 019.233SG.010.

\section*{Data availability}
The wildfire risk map~\cite{pyrologix}, ignition records~\cite{usfs_ignition_points}, and land-cover datasets~\cite{USCensusUrban2020,WFPI2020} used in this study are publicly available from the cited sources. A record of all datasets used and generated in this study is available at Zenodo under DOI \url{10.5281/zenodo.21359921}.

\section*{Code availability}
The infrastructure placement and drone routing models (Julia) and preprocessing/simulation scripts (Python), as well as the code and data to re-generate all the results and figures of this paper, are available at \danique{\url{https://github.com/RomainPuech/wildfire_drone_routing/tree/public-release}} under MIT license and archived at Zenodo under DOI \url{10.5281/zenodo.21359445}.

\raggedright
\bibliography{sn-bibliography}

\section*{Appendix A. Pyrologix predictive quality}
\label{app:pyrologix_quality}
We assess the static Pyrologix burn-probability map \cite{pyrologix} against California
ignition records spanning 2021 to 2024. The background median over all valid California
cells is $0.619$ on the 0--1 scale and is year-independent (the map is static).
Table~\ref{tab:pyrologix_quality} reports, for each year, the median Pyrologix value at
ignition cells, its ratio to the background median, and the share of ignitions whose cell
value exceeds that median. In 2021, $75.1\%$ of fires fall above the
background median, $+25.1$ percentage points above a uniform-random baseline. Across the
out-of-sample years 2022--2024, this share ranges from $73.7\%$ to $78.4\%$; no year
deviates more than $3.3$ percentage points from the 2021 value. The most recent
out-of-sample year, 2024, yields $73.7\%$ of fires above background, consistent with
the multi-year pattern. This stable predictive signal supports the use of the Pyrologix
surface as the risk layer for sensor placement and drone routing throughout the
experiments.

\begin{table}[!ht]
\caption{\textbf{Pyrologix predictive quality for California ignition locations 2021--2024.} The background median is computed over all valid California cells (Pyrologix is static so the background is year-independent). All values are on the 0--1 scale. For each year we report the median Pyrologix value at ignition cells, its ratio to the background, the share of ignitions whose cell value exceeds the background median, and the improvement over the random-baseline 50\%.}
\label{tab:pyrologix_quality}
{\small
\begin{tabular}{@{}lcccccc@{}}
\toprule
\textbf{Year} & \textbf{N fires} & \textbf{Fire median} & \textbf{Bg median} & \textbf{Ratio} & \textbf{\% above bg} & \textbf{Improvement (pp)} \\
\midrule
2021 & 981 & 0.683 & 0.619 & 1.104 & 75.1\% & 25.1 \\
2022 & 903 & 0.690 & 0.619 & 1.116 & 74.6\% & 24.6 \\
2023 & 992 & 0.702 & 0.619 & 1.134 & 78.4\% & 28.4 \\
2024 & 817 & 0.683 & 0.619 & 1.104 & 73.7\% & 23.7 \\
\bottomrule
\end{tabular}
}
\end{table}

\section*{Appendix B. Adaptations to the PSO routing heuristic}
\label{app:pso}

Our drone routing builds on the PSO-inspired heuristic for the Team Orienteering
Problem (TOP) of Dang, Guibadj and Moukrim~\cite{PSO}. We summarise that algorithm
briefly and then describe the changes we made so that it runs fast enough for
large, grid-based wildfire instances. The original method does not represent each
drone route separately. Instead, a candidate solution is encoded as a single
permutation of all the locations to be considered, called a giant tour. A
\emph{split} procedure reads this permutation in order and, using dynamic
programming, cuts it into the set of feasible routes (one per drone) that collects
the most reward while respecting the battery budget. The algorithm then improves
solutions in the style of particle swarm optimization, recombining promising
permutations and refining them with local search moves that relocate a single
location (\emph{shift}) or exchange two of them (\emph{swap}). Each candidate move
is scored by re-running the split, so the split and the local search together
account for almost all of the running time.

We instantiate this framework on the wildfire grid: the locations are grid cells
carrying a wildfire-risk reward, the depots are charging stations, and travel cost
is the number of grid steps (Chebyshev distance) between cells. We also extend the
formulation so that within each battery cycle a drone may depart from and return to
any charging station rather than a single fixed depot. The giant tour therefore
contains several depot markers, and the split is free to assign each route to a
different station.

The remaining changes are accelerations that exploit the grid structure without
altering which solutions the search explores. First, we observe that a route can
only begin at a charging-station marker, and these markers are few compared with
the number of grid cells. The original split nonetheless scans every position of
the permutation; we instead run the dynamic program over depot positions only,
which lowers the cost of evaluating one candidate from roughly $\mathcal{O}(n\,m)$
to about $\mathcal{O}(k\,m \log k)$, where $n$ is the number of cells, $m$ the
number of drones, and $k$ the (small) number of depot markers. For typical
instances this alone speeds up evaluation by around sevenfold. Second, during local
search we derive simple conditions under which a candidate swap or shift cannot
possibly change the split outcome (for example, moves confined to cells that no
route currently visits), and we skip the expensive re-evaluation of such moves.
Depending on how densely the routes fill the grid, this discards between $40\%$ and
$85\%$ of candidate moves with no effect on the result. Third, rather than rebuild
every route after each accepted move, we cache the routes and recompute only the
ones a move actually touches, skipping the dynamic program entirely when nothing
relevant changes (about three quarters of the time); this makes the local search
roughly four times faster. Finally, we replace per-lookup dictionary access to
travel costs with a dense distance matrix and reuse pre-allocated buffers in the
inner loops, removing memory-allocation overhead from the hot path.

Battery feasibility is enforced inside the split itself: a route stops collecting
cells as soon as the drone could no longer return to a charging station within its
budget, so the procedure never produces an infeasible route. A greedy step then
extends any route that the search left short, allowing each drone to use its full
battery budget before returning. Each of the accelerations above was verified to
leave solution quality unchanged relative to a direct re-evaluation, so the reported
routing results reflect the same search as the baseline method, obtained at
substantially lower computational cost. The exact implementation of the PSO-based routing heuristic is available in the code (see the Code Availability statement).

\end{document}